\documentclass[]{article}

\usepackage[a4paper,margin=3.5cm]{geometry}

\usepackage{graphicx}
\usepackage{subcaption}
\usepackage[utf8]{inputenc}
\usepackage{amsmath,amsfonts,amssymb}
\usepackage{bbm}
\usepackage{hyperref}
\usepackage{todonotes}
\usepackage{color}
\usepackage{tikz}
\usepackage{gnuplot-lua-tikz}
\usepackage{algorithm}
\usepackage{algorithmic}
\usepackage{booktabs}

\newcommand{\R}{\mathbbm{R}}
\newcommand{\Ft}{{\mathcal{F}_t}}
\newcommand{\LL}{{\mathcal{L}}}

\newcommand{\Gt}{{\Gamma_\text{top}}}
\newcommand{\Gi}{{\Gamma_\text{int}}}
\newcommand{\Gs}{{\Gamma_\text{side}}}
\newcommand{\Gb}{{\Gamma_\text{bottom}}}
\newcommand{\Oo}{ {\Omega_\text{out}} }
\newcommand{\Oone}{ {\Omega_{\text{cell}, 1}} }
\newcommand{\Oi}{ {\Omega_{\text{cell}, i}} }
\newcommand{\Oj}{ {\Omega_{\text{cell}, j}} }
\newcommand{\OiN}{ {\Omega_{\text{cell}, N}} }

\newcommand{\HOneZero}{{H^1_0(\Omega, \R^d)}}
\newcommand{\HOne}{{H^1(\Omega, \R^d)}}

\newcommand{\nuElas}{{\nu_\text{elast}}}
\newcommand{\nuVol}{{\nu_\text{vol}}}
\newcommand{\nuPeri}{{\nu_\text{peri}}}
\newcommand{\nuPenalty}{{\nu_\text{penalty}}}

\newcommand{\JElas}{{J_\text{elast}}}
\newcommand{\JVol}{{J_\text{vol}}}
\newcommand{\JPeri}{{J_\text{peri}}}

\providecommand{\keywords}[1]
{
  \small	
  \textbf{Keywords:} #1
}

\providecommand{\classification}[1]
{
  \small	
  \textbf{Mathematics Subject Classification (2010):} #1
}

\begin{document}

\title{A shape optimization algorithm for cellular composites}

\author{Martin Siebenborn\thanks{Universität Hamburg, Bundesstraße 55, 20146 Hamburg, martin.siebenborn@uni-hamburg.de} \and Andreas Vogel\thanks{Ruhr-Universität Bochum, Universitätsstraße 150, 44801 Bochum, a.vogel@rub.de}
}

\date{}

\maketitle

\begin{abstract}
We propose and investigate a mesh deformation technique for PDE constrained shape optimization. Introducing a gradient penalization to the inner product for linearized shape spaces, mesh degeneration can be prevented  within the optimization iteration allowing for the scalability of employed solvers. We illustrate the approach by a shape optimization for cellular composites with respect to linear elastic energy under tension. The influence of the gradient penalization is evaluated and the parallel scalability of the approach demonstrated employing a geometric multigrid solver on hierarchically distributed meshes.
\end{abstract}

\keywords{shape optimization, multigrid methods, algorithmic scalability\\[0.2cm]}
\classification{49Q10, 35Q93, 65M55, 65Y05}

\section{Introduction}
\label{sec:Introduction}

Partial differential equation (PDE) constrained shape optimization (e.g. \cite{Haslinger2003Introduction,delfour2001shapes,sokolowski1992introduction}) aims at finding optimal structural designs for models described by PDEs.
In contrast to classical optimal control, where model parameters are represented by distributed functions, the underlying assumption of shape optimization is that there is a fixed number of available parameters (e.g. materials) which are spatially separated by interfaces.
The variable in this optimization technique is then the shape of these contours while the topology is assumed to be known and kept fixed.

Nowadays, shape optimization has reached a high level of maturity, which can be seen in the field of applications ranging from, e.g., aerodynamics \cite{mohammadi2001applied,jameson2003aerodynamic,giles2000introduction} and acoustics \cite{berggren2016largescale} over electrostatics \cite{langer2015shape} to solid mechanics \cite{allaire2012shape}. 
The advantage over classical optimal control is a significant reduction of degrees of freedom allowing for high resolutions in regions of interest such as the interfaces between materials.
Moreover, the choice of function spaces for an optimal control typically introduces a certain level of smoothness to the parameter, which prevents desired sharp interface solutions.
The price one has to pay is to provide initial knowledge about the number of parameters, their approximate location, and topology.
Typically, shapes are numerically encoded by edges (2d) or faces (3d) in a finite element mesh.
This enables true integer solutions with sharp interfaces, but the initial topology can not be changed such that it is not straightforward to join or split shapes.

%Two alternative methodologies for these types of problems should be mentioned.
%Probably the most obvious one would be to introduce an integer variable, one for each finite element, which encodes material properties.
%Yet, solving this problem with methods from discrete optimization might significantly limit the spatial resolution.
%Moreover, one tends to obtain solutions, which are not stable under mesh refinements and possibly show checkerboard patterns.
For the special class of interface identification problems there are two prominent, alternative approaches to shape optimization:
first, interfaces can be implicitly resolved by level-set (cf.\ \cite{sethian1999level}) or phase-field functions (see for instance \cite{penzler_rumpf_wirth_2012}, where a similar problem to the one in this paper is studied).
Second, there is the topology optimization approach, where the binary variable describing material and void is relaxed and a combination of continuous optimization and filter techniques are applied \cite{bendsoe2011topology}.
Similarly to our considerations in this paper, there is also ongoing effort in maintaining material thicknesses in the context of level-set methods (see for instance \cite{allaire2016thickness}).

The focus of this article are applications where the shape topology should intentionally be retained. A representative example is a biological cell composite where cells may deform and move relative to each other, but the topology of cells remains unchanged. E.g., the cell layout of the epidermis (the outermost layer of the human skin) consists of different (and differently shaped) cells surrounded by thin lipid layers. We refer the reader to, e.g., \cite{ADDRSkin2013,Querleux2014,Naegel09,Limbert2017MathematicalAC} for a review of modeling aspects. For such a system, a reasonable question is to ask for an optimized cell arrangement with respect to elastic energy under outer force application. Using shape optimization to analyze such objectives, a certain numerical difficulty is due to the thin layers between cells that must be treated appropriately and we develop mathematical tools to address such settings.  

%Here we study methodologies towards elasiticity models for the outermost layers of the human skin.
%Thus, we distinguish two different materials with different elastic behavior, i.e. skin cells and the lipid layers separating them.

In order to optimize shapes it is obligatory to define a set of admissible shapes and a corresponding tangent space for the descent directions.
The probably most prominent approach is to describe shapes as points on a Riemannian manifold \cite{michor2005vanishing,michor2006riemannian}.
A lot of attention has been paid to finite element approximations of this manifold and the corresponding tangent spaces (see for instance \cite{langer2015shape,schulz2016computational,schulz2016efficient,iglesias2018two}).

In this paper, we study a modification of inner products in the linearized FEM approximation of shape spaces.
We propose a formulation that adapts the modeling of biological tissues to some extend.
During shape updates in the optimization process, descent direction with small gradients should be promoted, since
in growth processes of biological cellular structures the material can only be compressed to some extend and not completely vanish. This can be modeled by limiting gradients of shape updates.
A resulting feature of limited deformation gradients is the preservation of the FEM mesh quality, which is demonstrated within this work. In previous work \cite{siebenborn2017algorithmic} it has been observed that the mesh quality may reduce dramatically during shape updates with the consequence of unreliable finite element approximations.
We also demonstrate that the inner product in shape space can be implemented using standard finite element techniques, which is our foundation for scalable algorithms on parallel supercomputers. 

The main contributions of this work thus are:
\begin{itemize}
  \item Optimization algorithm with gradient penalization for cellular composites
  \item Numerical study of the influence of the gradient penalization on optimization results, mesh quality, and solver behavior for a cell composite optimization
  \item Parallel scalability analysis for the presented implementation
\end{itemize}

This paper is structured as follows:
In Section~\ref{sec:problem} the PDE model is introduced and the shape optimization problem is formulated.
Section~\ref{sec:optimization-algorithm} describes a gradient based optimization algorithm and the novel inner product is proposed.
Numerical experiments and parallel performance of our implementation are then investigated in Section~\ref{sec:Numerics}.
Finally, Section~\ref{sec:conclusion} discusses the achieved results and provides a brief outlook on future research.

\section{Problem formulation}
\label{sec:problem}

\begin{figure}
\begin{center}
\def\svgwidth{0.6\textwidth}
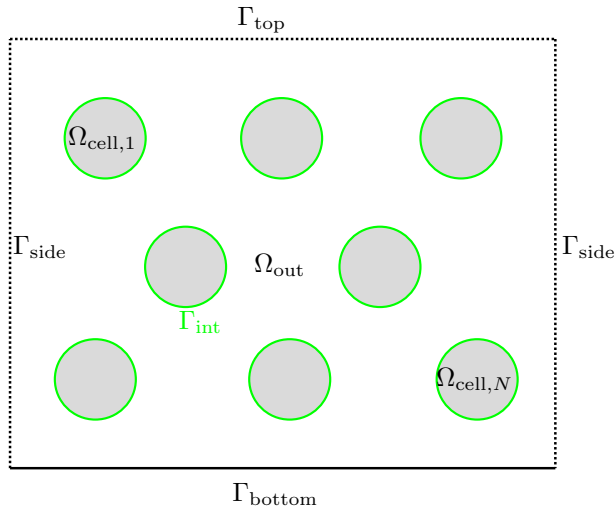
\end{center}
\caption{Schematic view of the domain $\Omega$ consisting of an bulk material $\Oo$ and cell inclusions $\Oone, \dots, \OiN$ with different material properties}
\label{fig:domain}
\end{figure}

The general idea of our model can be illustrated by a geometrical setting as sketched in Figure~\ref{fig:domain}: we consider a composite material that is fixed at its bottom surface and forces act on the top surface. We assume that the domain of interest consists of several material inclusions with different elastic behavior. This setting is inspired by a composite of biological cells with some surrounding bulk material. The goal is then to improve the stiffness of the entire object by changing its interior geometrical configuration, i.e. relocating and reshaping the inclusions. As an important constraint, the outer shape of the block has to remain unchanged in order to prevent degenerated solutions such as removing the force by a vanishing upper surface. If the bulk material and each of the inclusions consist of a homogeneous material, the stiffness of the composite block is a function of the shapes of the inclusions, and our aim is to devise an algorithm for the efficient computation of an optimized shape of these cells. In the following, we provide a precise mathematical formulation for this setting.

\subsection{Model equations}
\label{sec:model}

Let $\Omega \subset \R^d$ and $\Oo, \Oone, \dots, \OiN \subset \Omega$ ($N\in \mathbbm{N}$) be bounded Lipschitz domains in a $d$-dimensional space for $d \in \lbrace 2, 3\rbrace$.
We assume that the entire domain is made up by the bulk material $\Oo$ and the cell inclusions $\Oi$, $1 \leq i \leq N$, and that cell inclusion do not touch each other, i.e. we assume 
\begin{align}
\Oo = \Omega \setminus \left(\bigcup\limits_{i=1}^{N} \mathrm{cl}(\Oi)\right)
\quad \text{ and } \quad
\mathrm{cl}(\Oi) \cap \mathrm{cl}(\Oj) = \emptyset \text{ for } i \neq j,
\label{eq:domain-definition}
\end{align}
where $\mathrm{cl}(\Oi)$ denotes the closure of a set.
This implies that there is always material between two cells $\Oi$ and $\Oj$. 
The boundary of the cells is denoted by $\Gi := \cup_{i=1}^N \partial \Oi$.
% for the inclusions $\mathrm{cl}(\Oi)$ are pairwise disjoint for $i=1,\dots,N$ and that $\Omega = \Oo \cup \left(\bigcup\limits_{i=1}^{N} \mathrm{cl}(\Oi)\right)$.
The outer boundary of $\Omega$ is subdivided into the disjoint parts $\Gt, \Gb, \Gs$, each with a positive Lebesgue measure.
An example of this situation is depicted in Figure \ref{fig:domain}.

The mathematical model is given by the equations of linear elasticity for a displacement field $u: \R^d \to \R^d$ as
\begin{equation}
\begin{aligned}
\mathrm{div}(\sigma(u)) = 0, \quad &\text{in}\; \Omega,\\
\sigma(u) \cdot n = f, \quad &\text{on}\; \Gt,\\
\sigma(u) \cdot n = 0, \quad &\text{on}\; \Gs,\\
u = 0, \quad &\text{on}\; \Gb,
\end{aligned}
\label{eq:elasticity-strong}
\end{equation}
where $\sigma(u) = \lambda \, \mathrm{tr}(\nabla u) I + \mu (\nabla u + \nabla u^T)$ is the stress tensor in terms of the Lamé parameter $\lambda$ and $\mu$, and $f: \R^d \to \R^d$ is a surface force.

Defining a function space $\mathcal{V} := \left\lbrace u \in \HOne:\, u = 0 \; \text{a.e. on} \; \Gb \right\rbrace$ on $\Omega$ we obtain the weak formulation of problem \eqref{eq:elasticity-strong} as follows:
Find $u \in \mathcal{V}$ such that
\begin{equation}
\int_\Omega \sigma(u) : \nabla w \, dx = \int_\Gt f \cdot w\, ds \quad \forall\, w \in \mathcal{V}.
\label{eq:elasticity-weak}
\end{equation}

The purpose of this model is the simulation of the elastic behavior of a composite material.
We assume that $\Oo$ is homogeneously filled with one material.
In addition, each of the $N$ inclusions $\OiN$ is filled with a homogeneous material, but the material may differ between inclusions. This situation is reflected in the mathematical model by piece-wise constant Lamé parameter $\lambda, \mu \in L^\infty(\Omega)$.

\subsection{Elasticity optimization problem}
\label{sec:optimization-problem}

For the above settings, we model the stiffness of the material as a function of the inner geometrical configuration.
We thus consider a minimization problem as 
\begin{equation}
\begin{aligned}
\min\limits_{\Gi}\quad & \JElas(u, \Omega) :=  \frac{1}{2}\int_\Omega \sigma(u):\nabla u\, dx\\
\text{s.t.}\quad &\int_\Omega \sigma(u) : \nabla w \, dx = \int_\Gt f \cdot w\, ds \quad \forall\, w \in \mathcal{V}.
\end{aligned}
\label{eq:optimization-problem}
\end{equation}
This minimization problem is to be understood in a formal way, i.e.
we optimize over the interior shape $\Gi$ of the domain $\Omega$. Note that a change of the location and shape of the inclusions $\Oone, \dots, \OiN$ changes the spatial distribution of the Lamé parameter and thereby affects the objective function in \eqref{eq:optimization-problem}.

In order to represent shape modifications we concentrate on the so-called perturbation of identity.
This means that, in terms of the initial shape $\Omega = \Omega_0$, we considered for $t\geq 0$ deformed domains as
\begin{equation}
\Omega_t = \lbrace \Ft(x): \; x \in \Omega \rbrace
%\end{equation}
\quad \text{ with } \quad
%\begin{equation}
\begin{aligned}
\Ft: \Omega \to \Omega_t,\;
\Ft (x) = x + tv(x)
\end{aligned}
 \label{eq:perturbation-of-identity}
\end{equation}
for a given displacement vector field $v \in W^{1, \infty}_0(\Omega, \R^d)$. 
Note that this choice guarantees that the outer shape of $\Omega$ is retained.
Moreover, for $t$ sufficiently small, $\Omega_t$ remains admissible in the sense that there are no overlappings of the cell inclusions.

In order to investigate sensitivities of the objective $J$ with respect to the shape we mainly follow the concepts introduced in \cite{sokolowski1992introduction} and \cite{delfour2001shapes}.
For a vector field $v$ as defined above, the directional shape derivative of the functional $J$ evaluated at $\Omega$ in direction $v$ is defined by
\begin{equation}
dJ(\Omega)[v] = \lim\limits_{t \searrow 0} \frac{ J(\Omega_t) - J(\Omega)}{t}.
\label{eq:directional-derivative}
\end{equation}
Note that the general definition above implicitly takes into account the dependence of $J$ on the state variable $u$.
Since $u$ is the solution of the elliptic equation \eqref{eq:elasticity-weak} it depends on the interior shape of the domain $\Omega_t$, in particular on the distribution of Lamé parameters.
Therefore, $J(u_t, \Omega_t)$ has to be considered in \eqref{eq:directional-derivative} and derivatives with respect to $t$ have to be computed.

We consider the associated Lagrangian to problem \eqref{eq:optimization-problem} given by
\begin{equation}
\LL(u,p, \Omega) = \frac{1}{2} \int_\Omega \sigma(u):\nabla u\, dx - \int_\Omega \sigma(u):\nabla p \, dx + \int_\Gt f p \, dx.
\label{eq:lagrangian}
\end{equation}
The state variable $u$ and the corresponding adjoint $p$ can then be characterized as the saddle point of \eqref{eq:lagrangian} as
\begin{equation}
J(u, \Omega) = \min\limits_{u \in \HOne} \max\limits_{p \in \HOne} \LL(u,p,\Omega).
\label{eq:min-max}
\end{equation}
In order to obtain the shape derivative $dJ(\Omega)[v]$, the equation above is reformulated in terms of $\Omega_t$, $u_t$ and $p_t$. 
Following the presentations in \cite[Chapter 10]{delfour2001shapes}, where an application of the Correa and Seeger theorem \cite[Theorem 2.1]{correa1985directional} to shape sensitivity analysis is described, the differentiability of the right hand side in \eqref{eq:min-max} with respect to $t$ is obtained.
From these results it follows that 
\begin{equation}
dJ(\Omega)[v] = \frac{d}{dt} \left.\LL (u_t, p_t, \Omega_t)\right|_{t=0}.
\end{equation}

First we observe that in a saddle point for $t=0$ it holds $\LL_u(u,p,\Omega) = 0$.
This yields the condition $p = u$ for the adjoint variable, which we can therefore substitute in the following.
Let $u^t := u_t \circ \Ft$ be the pull-back of $u_t$ onto the domain $\Omega$ and observe that $\nabla u_t = (D\Ft)^{-T} \nabla u^t$. 
In addition, due to the special choice of $\Ft$ it holds
\begin{equation}
\begin{aligned}
D\Ft = I + t\nabla v, \quad %\\
\left.\frac{d}{dt} \det(D \Ft)\right|_{t=0} = \mathrm{div}(v), \quad %\\
\left.\frac{d}{dt} (D \Ft)^{-1} \right|_{t=0} = -Dv.
\end{aligned}
\end{equation}
With these identities we can compute
\begin{multline*}
 \frac{d}{dt}\left.\LL (u_t, p_t, \Omega_t)\right|_{t=0}
= \frac{d}{dt}\left. \frac{1}{2} \int_{\Omega_t} \sigma(u_t):\nabla u_t \, dx \right|_{t=0}\\
= \frac{d}{dt}\left. \frac{1}{2} \int_{\Omega_t} \left( \left[(\lambda\circ \Ft^{-1})\, \mathrm{tr}(\nabla u_t)I + (\mu \circ \Ft^{-1})\left( \nabla u_t + \nabla u_t^T\right)\right]:\nabla u_t \right)  \, dx \right|_{t=0}\\
= \frac{d}{dt} \frac{1}{2} \int_{\Omega} \left[\lambda\, \mathrm{tr}((D\Ft)^{-T} \nabla u^t)I + \mu\left( (D\Ft)^{-T}\nabla u^t + ((D\Ft)^{-T}\nabla u^t)^T\right)\right]\\
  \left.:((D\Ft)^{-T}\nabla u^t)\, \det(D\Ft) \, dx \right|_{t=0}\\
= \frac{1}{2}\int_{\Omega} \left[\lambda\, \mathrm{tr}((-Dv)^T \nabla u)I + \mu\left( (-Dv)^T\nabla u + ((-Dv)^T \nabla u)^T\right)\right]:\nabla u \\
+ \sigma(u):\left( (-Dv)^T \nabla u \right) + \mathrm{div}(v) \left(\sigma(u):\nabla u \right) \, dx
\end{multline*}
We finally obtain the directional shape derivative
\begin{align}
d\JElas(u, \Omega)[v] 
%:= dJ(u, \Omega)[v] 
&= \int_\Omega \frac{1}{2} \mathrm{div}(v) \left( \sigma(u):\nabla u \right) - \sigma(u) : \left(\nabla v \nabla u \right)\, dx \notag \\
&= \int_\Omega  \left( \frac{1}{2} \left(\sigma(u):\nabla u\right) I -  (\nabla u)^T \sigma(u) \right) : \nabla v \, dx,
\label{eq:shape-deriv-elasticity}
\end{align}
where we exploited the special structure of the stress tensor for isotropic linear elastic material. 
This expression can also be found in \cite{novotny2007topological} with a slightly different derivation. 
Note that this formula is different to the expression used in the previous work \cite{siebenborn2017algorithmic} where we found numerical artifacts influencing the optimization. In this paper, we use the above expression and modify the algorithm in order to address these issues.

%An immediate consequence of the sensitivity analysis applied to the model problem introduced in \eqref{eq:optimization-problem} would be an unintended behavior close to the boundaries and particularly $\Gt$ and $\Gb$.
Maintaining the outer shape of $\Omega$ is not yet incorporated into the model.
Quite the contrary: the sensitivities $dJ(\Omega)[v]$ might show largest absolute values close to the boundaries due to the fact that shrinking $\Gt$ or extending $\Gb$ typically leads to a steep descent of the objective. To avoid this unintended behavior, we thus define a domain $D \subset \Omega$, slightly smaller than $\Omega$, which leaves out the outer boundaries.
The directional derivatives $dJ(\Omega)[v]$ are then only computed within $D$.
Note that this procedure is equivalent to considering deformation fields $v$ with support in $D$ only.
In addition, from a practical aspect, one typically encounters nonzero values in the discrete approximation of $dJ(\Omega)[v]$ even for finite element shape functions that do not intersect $\Gi$.
However, a FE shape function $\phi: D \to \R^d$ can not contribute to a descent direction $v$ if $\mathrm{supp}(\phi) \cap \Gi = \emptyset$.
%to be one FE shape function used for the discretization of $v$.
%In Section~\ref{sec:Numerics}, e.g., $\phi$ is a piece-wise linear functions on a patch of triangles.
%Observe that, if $\mathrm{supp}(\phi) \cap \Gi = \emptyset$, this shape function can not contribute to a descent direction $v$.
Since Lamé parameter are constant for each cell and the bulk material, only deformations of $\Gi$ can affect the objective functional. We thus add an additional step to the optimization algorithm (cf. Alg.~\ref{alg:optimization-algorithm}), which sets these shape derivative values to zero in order to avoid unintended deformations. 

\subsection{Optimization with respect to geometrical aspects}
\label{sec:geometric-optimization-problem}

Our problem setting is motivated by an elasticity model for biological cell composites such as the human skin. For such a system, also geometrical optimization aspects play a role. 
We therefore modify the objective function by two weighted additional terms using $\nuElas >0$ and $\nuVol, \nuPeri \geq 0$, and defining
\begin{equation}
\begin{aligned}
\min\limits_{\Gi}\quad & J(u, \Omega) := \nuElas \JElas(u, \Omega) + \nuVol \JVol(\Omega) + \nuPeri \JPeri(\Omega)\\
& \quad\quad\quad\;\;=  \frac{\nuElas}{2} \int_\Omega \sigma(u):\nabla u\, dx + \nuVol \int_{\Oo}1\, dx +  \nuPeri \int_{\Gi} 1\, ds\\
\text{s.t.}\quad &\int_\Omega \sigma(u) : \nabla w \, dx = \int_\Gt f \cdot w\, ds \quad \forall\, w \in \mathcal{V}.
\end{aligned}
\label{eq:geometric-optimization-problem}
\end{equation}
The additional terms represent the minimization of the volume of the bulk material in $\Oo$ and the minimization of the perimeter (2d) or surface (3d) of the inclusions $\Gi$.
The motivation for these terms is to find a space-filling design with minimal cell perimeter/surfaces. %depending on the parameters $\nuElas, \nuVol, \nuPeri$.
For a derivation of the shape derivatives in direction $v$ of the geometrical quantities, given by
\begin{align}
d\JVol(\Omega)[v] &= \int_{\Gi} v\cdot n\, ds = - \int_{\Oo} \mathrm{div}(v) \,dx, 
\label{eq:shape-deriv-volume}
\\
%\qquad 
d\JPeri(\Omega)[v] &= \int_{\Gi} \mathrm{div}(v) - \frac{\partial v}{\partial n} \cdot n\, ds,
\label{eq:shape-deriv-perimenter}
\end{align}
we refer the reader to \cite{sokolowski1992introduction}.
In the equation above, $n$ is assumed to be the normal vector field at $\Gi$, which points outwards of the inclusions $\Oone, \dots, \OiN$. 

\section{Optimization algorithm}
\label{sec:optimization-algorithm}

The equations~\eqref{eq:shape-deriv-elasticity}, \eqref{eq:shape-deriv-volume}, and \eqref{eq:shape-deriv-perimenter} state the sensitivities of the objective in terms of variations according to the perturbation of identity.
Since our aim is a gradient descent method, a space of admissible shapes and an according inner product has to be determined. The shape derivatives can then be represented as gradients with respect to this inner product and used as shape updates in the sense of a descent direction.

It is a common approach (see, e.g., \cite{langer2015shape,schulz2016efficient}) to locally approximate the shape space by $\HOneZero$-deformations to the initial shape $\Omega$ and equip this with the bilinear form originating from an elliptic PDE.
Here, we focus on linear elasticity, i.e.\ we choose an inner product on $\HOneZero$ by
\begin{equation}
a(v,w) := \int_\Omega \sigma(v):\nabla w \, dx.
\end{equation}

Let $dJ(\Omega)[w]$ be the shape derivative evaluated in $\Omega$ in direction $w$.
In order to obtain the deformation field $v \in \HOneZero$ we then solve
\begin{equation}
a(v,w) = dJ(\Omega)[w] \quad \forall \; w \in \HOneZero.
\label{eq:inner-product}
\end{equation}
This step is the representation of the derivative $dJ(\Omega)[\cdot]$ with respect to the inner product $a(\cdot, \cdot)$ in $\HOneZero$.
Note that deformations in $\HOneZero$ are not admissible in general, since a Lipschitz domain deformed via $\Ft$ in terms of $v$ may loose this property.

A consequence of a gradient descent method applied to a situation as the one depicted in Figure \ref{fig:domain} is that overlappings are likely to appear during the optimization.
Moreover, since the domain $\Omega$ and the inclusions $\Oone, \dots, \OiN$ are represented by edges and faces in a finite element mesh, the topology can not be changed.
It is thus of interest to modify $a(\cdot, \cdot)$ such that large gradients in the resulting deformation field are avoided.
This has two effects: first, the quality of the FEM mesh is retained to some extend even under large deformations. This is crucial for the applied linear solver and the overall approximation quality of the FEM method.
Second, we can move modeling aspects into $a(\cdot, \cdot)$.
The ultimate goal of our method is to reflect the behavior of biological tissue in the human skin.
Selecting descent directions which avoid large gradients promotes shape updates where the thickness of channels between cells does not tend to zero too fast.

Recall that \eqref{eq:inner-product} is the necessary and sufficient optimality condition of the minimization problem

\begin{equation}
\min\limits_{v \in \HOneZero} j_1(v) 
:= \frac{1}{2} a(v,v) - dJ(\Omega)[v].
%:= \frac{1}{2}\int_\Omega \sigma(v):\nabla v \, dx - dJ(\Omega)[v].
\label{eq:elastic-minimization}
\end{equation}
Our aim is to modify the selected descent directions to prefer translations over compression of the material in $\Omega$.
We thus enrich \eqref{eq:elastic-minimization} by a penalty term involving the gradients of the shape update.
We thus propose the following modified formulation
\begin{equation}
\min\limits_{v \in \HOneZero} j(v) 
:= j_1(v) + j_2(v) 
:=
 j_1(v) 
+ \frac{\nuPenalty}{4} \int_\Omega \left( \left( \nabla v:\nabla v - b^2 \right)^+ \right)^2\, dx
% + \frac{\nuPenalty}{2} \left\Vert \left( \frac{1}{2} \Vert \nabla v \Vert^2_F - \frac{1}{2} b^2 \right)^+ \right\Vert_{L^2(\Omega)}^2
\label{eq:elastic-gradient-minimization}
\end{equation}
for a bound $b>0$ on the Frobenius-norm of the gradient and a penalty factor $\nuPenalty>0$.
Here, $(\cdot)^+$ denotes the positive part, i.e.\ $\max \lbrace 0, \Vert \nabla v \Vert^2_F - b^2 \rbrace$ with $\Vert \nabla v \Vert_F^2 = \nabla v: \nabla v$.
For large penalty parameter $\nuPenalty$, gradients exceeding a certain threshold are thus avoided.
%Observe that 
%\begin{equation}
%\frac{\nuPenalty}{2} \left\Vert \left( \frac{1}{2} \Vert \nabla v \Vert^2_F - \frac{1}{2} b^2 \right)^+ \right\Vert_{L^2(\Omega)}^2 =
%\frac{\nuPenalty}{2} \int_\Omega \left( \left( \frac{1}{2} \nabla v: \nabla v - \frac{1}{2} b^2 \right)^+ \right)^2\, dx
%\end{equation}
% \begin{equation}
% j_2(v) =
% \frac{\nuPenalty}{2} \int_\Omega \left( \left( \frac{1}{2} \nabla v: \nabla v - \frac{1}{2} b^2 \right)^+ \right)^2\, dx
% \end{equation}
%from which 
A generalized derivative evaluated in a direction $w\in \HOneZero$ follows as
%\begin{multline}
%\frac{d}{d v} \left[ \frac{\nuPenalty}{2} \left\Vert \left( \frac{1}{2} \Vert \nabla v \Vert^2_F - \frac{1}{2} b^2 \right)^+ \right\Vert_{L^2(\Omega)}^2 \right](w)\\
%= \nuPenalty \int_\Omega \left( \frac{1}{2} \nabla v: \nabla v - \frac{1}{2} b^2 \right)^+ \, \nabla v : \nabla w \, dx.
%\end{multline}
\begin{equation}
dj_2(v)[w]
= \nuPenalty \int_\Omega \left( \nabla v: \nabla v - b^2 \right)^+ \, \nabla v : \nabla w \, dx.
\end{equation}
We can then again formulate the necessary optimality conditions for \eqref{eq:elastic-gradient-minimization} and obtain the variational formulation:
Find $v\in \HOneZero$ such that
% \begin{multline}
% j^{\prime}(v)w =\\
% \int_\Omega \sigma(v):\nabla w \, dx - dJ(\Omega)[w] + 2 \nuPenalty \int_\Omega \left( \nabla v: \nabla v - b^2 \right)^+ \, \nabla v : \nabla w \, dx = 0
% \label{eq:nonlinear-PDE}
% \end{multline}
% for all $w \in \HOneZero$.
% For a better comparison to the inner product in \eqref{eq:inner-product}, the deformation field $v$ for one optimization step is obtained as the solution of the nonlinear variational formulation
\begin{equation}
g(v,w) := a(v,w) + 
dj_2(v)[w]
% 2 \nuPenalty \int_\Omega \left( \nabla v: \nabla v - b^2 \right)^+ \, \nabla v : \nabla w \, dx
 = dJ(\Omega)[w]
\quad \text{ for all } w \in \HOneZero.
 \label{eq:modified-inner-product}
\end{equation}
Note that $g(\cdot,\cdot)$ is nonlinear and not even differentiable.
The left hand side $g(\cdot,\cdot)$ only constitutes a suitable inner product in $\HOneZero$, if the term involving the positive part is inactive, i.e.\ when $\nabla v: \nabla v - b^2 \leq 0$ almost everywhere.
If gradients of $v$ exceed the threshold $b$, the formulation \eqref{eq:modified-inner-product} switches to a nonlinear behavior and can no longer be represented by a symmetric, coercive bilinear form  such as $a(\cdot, \cdot)$.% like equation \eqref{eq:inner-product}.

\begin{algorithm}[tb]
\caption{Gradient-penalized optimization algorithm}
\label{alg:optimization-algorithm}
\begin{algorithmic}[1]
\STATE Choose initial domain $\Omega^0$; choose $\nuElas, \nuVol, \nuPeri$; choose $\nuPenalty, b, t$
%\STATE $k \gets 0$
\FOR{$k = 0, 1, 2, \dots $}
  \STATE Compute displacement field $u^k$, i.e. solve model equations \eqref{eq:elasticity-weak}, for the current $\Omega^k$
  \STATE Compute shape derivative as given by the equations~\eqref{eq:shape-deriv-elasticity}, \eqref{eq:shape-deriv-volume}, \eqref{eq:shape-deriv-perimenter}, i.e. assemble
  $dJ(\Omega^k) = \nuElas d\JElas(u^k, \Omega^k) + \nuVol d\JVol(\Omega^k) + \nuPeri d\JPeri(\Omega^k)$%by eq.~\eqref{eq:geometric-optimization-problem}
  \STATE Reset shape derivative away from surface, i.e. set $dJ(\Omega^k)[\phi] = 0$ if $\mathrm{supp}(\phi) \cap \Gi = \emptyset$
  \STATE Compute descent direction $v$, i.e. solve $g(v,w) = dJ(\Omega^k)[w]$ for all $w \in \HOneZero$ 
  \STATE Update shape: $\Omega^{k+1} = \lbrace x + tv(x):\, x \in \Omega^k \rbrace$
% \STATE $k \gets k+1$
\ENDFOR
\end{algorithmic}
\end{algorithm}

Finally, we solve \eqref{eq:modified-inner-product} for $v$ via a semi-smooth Newton method and use the derivative in directions $w,h \in \HOneZero$ given by
\begin{multline}
j^{\prime \prime}(v)[w,h] = \int_\Omega \sigma(h):\nabla w \, dx\\
  + \nuPenalty \int_\Omega \gamma \left[ 2 \left(\nabla v : \nabla h \right) \cdot \left( \nabla v : \nabla w \right ) + \left( \nabla v: \nabla v - b^2 \right) \, \nabla h : \nabla w \right]\, dx,
\end{multline}
where 
\begin{equation}
L^\infty(\Omega) \ni \gamma =
\begin{cases}
1,\;\text{if}\;\Vert \nabla v \Vert_F^2 - b^2 > 0\\
0,\;\text{else}.
\end{cases}
% NOTE: I know that the square is not mandatory, since without it is equivalent. But it increases readablity
\end{equation}
%Parameters are chosen to be $\nuPenalty_0 = 1\mathrm{e}-1, \nuPenalty_\text{inc} = 1\mathrm{e}+1, \nuPenalty_\text{max} = 1\mathrm{e}+7, \lambda_\text{int} = 0.2, \mu_\text{int} = 1.5, \lambda_\text{out} = 0.02, \mu_\text{out} = 0.1, b = 1\mathrm{e}-1, \nu_\text{vol} = 1\mathrm{e}+1, \nu_\text{elas} = 1\mathrm{e}+2$.
%In the deformation equation we choose $\lambda = 0.1, \mu = 1.0$.

For theoretical background on Newton's methods for semi-smooth functions we refer the reader to \cite{ulbrich2002semismooth}.
Furthermore, the threshold on gradients of state variables is inspired by the theory of optimal control with gradient constraints (see for instance \cite{hintermller2010constrained}). Algorithm~\ref{alg:optimization-algorithm} shows the resulting overall algorithm. 

Since the outer shape of the domain $\Omega$ should not be modified, we restricted ourselves to a $\HOneZero$-approximation of descent directions in the shape space. 
However, for the numerical results in the next section, we slightly change the space for the gradient representation: we want to allow tangential movement of nodes along the boundary to allow for better mesh quality. % and more realism in the experiment.
In the two dimensional setting of Figure \ref{fig:domain} it is thus only required that the $x$-component at $\Gs$ and the $y$-component at $\Gb \cup \Gt$ is zero.
Note that this is sufficient to guarantee an unchanged outer boundary $\partial \Omega$.
% but allows for a broader behavior in the shape updates.
%In \eqref{eq:modified-inner-product}, the inner product is defined for deformation that completely vanish at $\partial \Omega$.

\section{Numerical Results}
\label{sec:Numerics}

We test the proposed method using a setup that is inspired by the topology of the upper human skin (see e.g. \cite{Naegel09} and the references therein). We emphasize that the presented numerical results are far from being biologically meaningful. In contrast, these simulations are solely used to demonstrate the numerical properties of the method, and we have chosen artificial material parameters and domain sizes for this example. However, the setup exposes a reasonable similarity with biologically observed patterns and we leave more realistic simulations for future work. 

\begin{figure}[tb]
  \centering
  \includegraphics[width=0.5\textwidth]{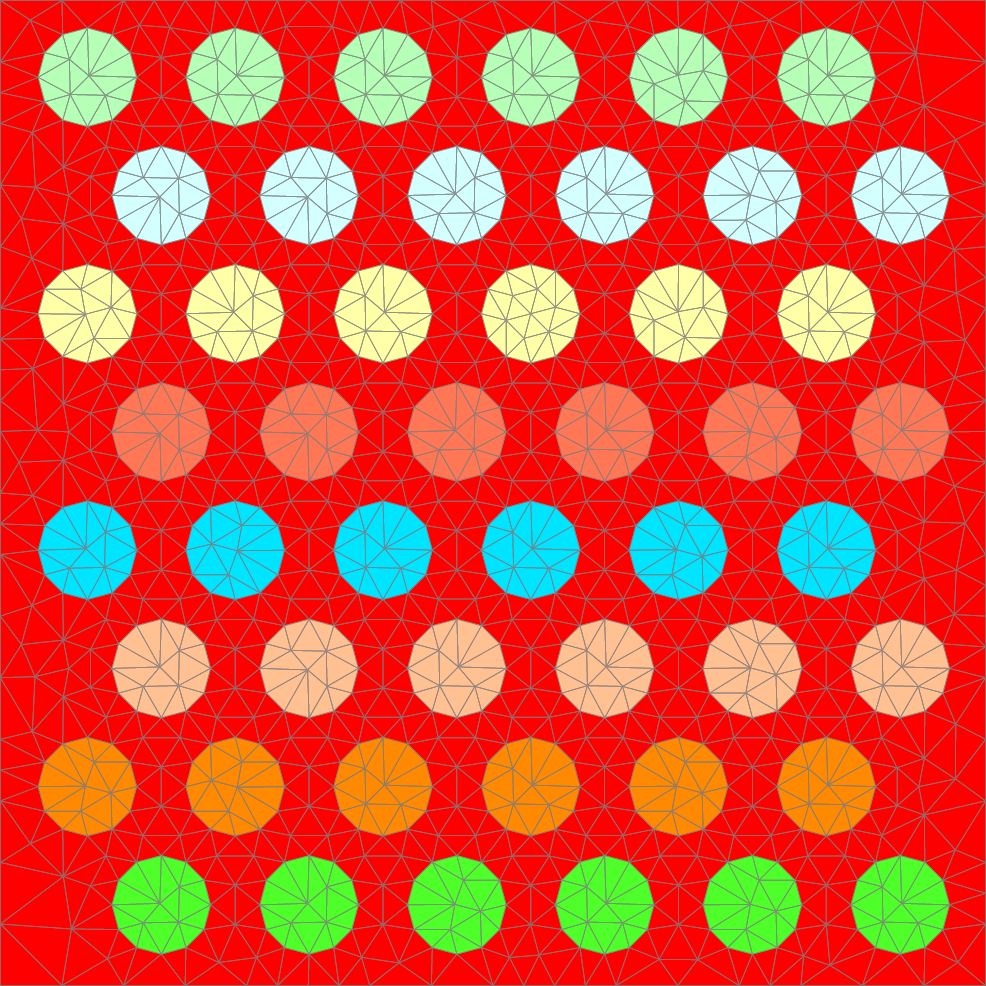}
  \caption{Initial coarse mesh. Different colors show different material properties}
  \label{fig:coarse-mesh-topology}
\end{figure}

For our test case, the topology of the domain $\Omega = (0,1)^2$ is shown in the initial configuration in Figure~\ref{fig:coarse-mesh-topology}. The outer domain (red) contains several cell inclusions which differ in their material properties. Cells of the same color have the same material properties and we number the eight cell types from top to bottom ($k=0,\dots,7$). For our numerical example, we choose the following material parameters
\begin{align}
  (\lambda, \mu) = (1.0, 0.1), \qquad & \text{ on } \Omega_\text{out},
\end{align}
for the outer space material, and the different cell types $\Omega_\text{cell,k}, k = 0,\dots,7$ are assigned a constant material parameter per cell, but with a linear increase in stiffness from top to bottom, given by
\begin{align}
 & (\lambda, \mu) = (\lambda_k, \mu_k), \text{ on } \Omega_\text{cell,k}, \text{ with } 
 \\
 \lambda_k &= (1 - \frac{k}{7}) \cdot  \lambda_\text{top} + \frac{k}{7}  \cdot \lambda_\text{bottom} , &&
  \lambda_\text{top} = 1.2, 
  \lambda_\text{bottom} = 2.0
 \\
 \mu_k &= (1 - \frac{k}{7}) \cdot \mu_\text{top} +  \frac{k}{7} \cdot \mu_\text{bottom}, &&
 \mu_\text{top} = 0.12, \mu_\text{bottom} = 0.2
\end{align}
For the elasticity problem, we apply a pulling force, $f = 0.1$, at the top boundary $\Gamma_\text{top}$, and use the plain strain assumption. For the computation of the shape derivate, we use globally constant values $(\lambda, \mu) = (0.1, 1.0)$ on the whole domain $\Omega$.
The coarse grid meshing has been performed with ProMesh \cite{ProMesh}, visualization is by ParaView \cite{Paraview}, and all simulations are carried out with the simulation toolbox UG4 \cite{Vogel2014ug4}.

\subsection{Qualitative optimization results}

\begin{figure}[tbp]
    \centering
    \begin{subfigure}[t]{\textwidth}
        \centering
        \makebox[0.04\textwidth][l]{(a)}
        \includegraphics[width=0.25\textwidth]{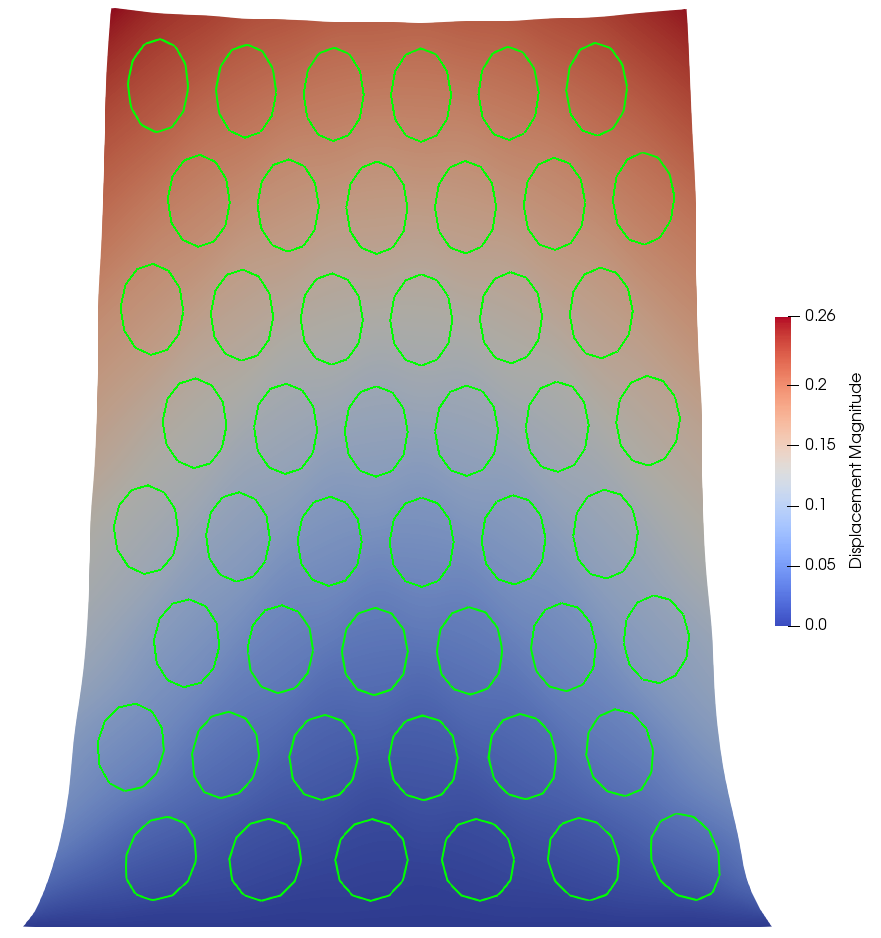}
        %\hfill
        \hspace*{1cm}
        \includegraphics[width=0.25\textwidth]{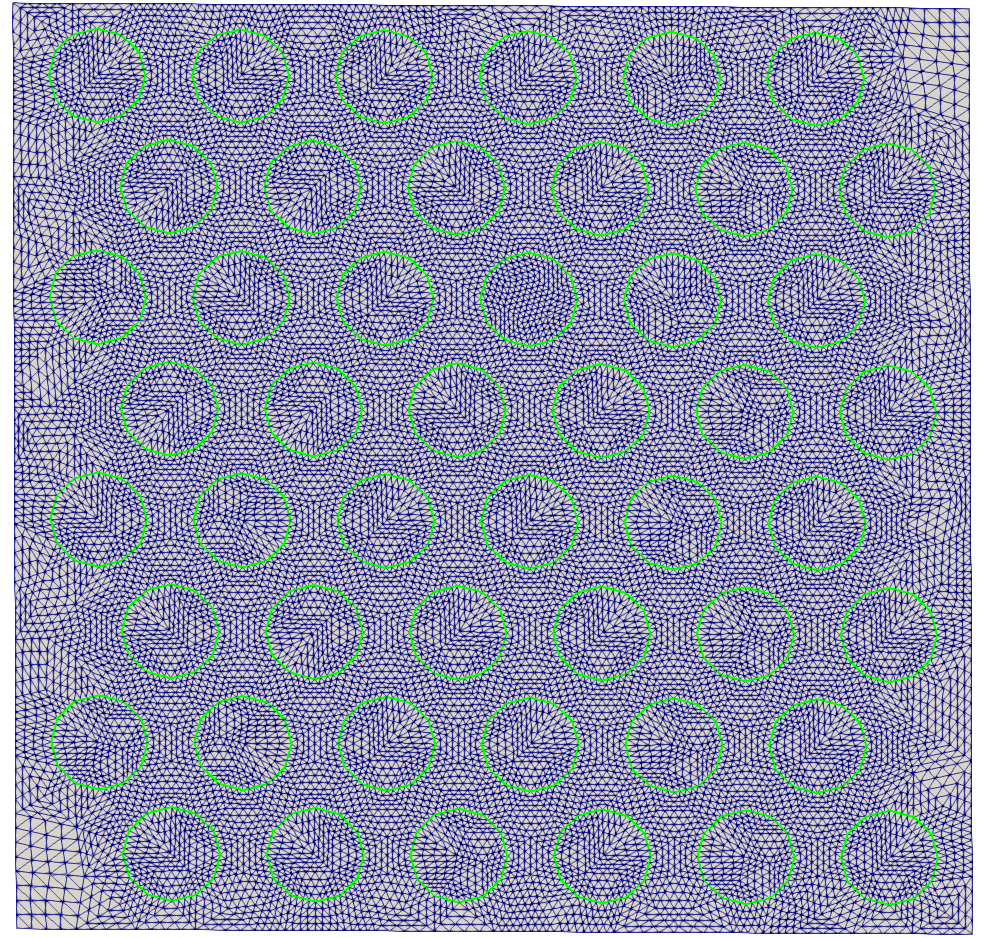}
        \includegraphics[width=0.2\textwidth]{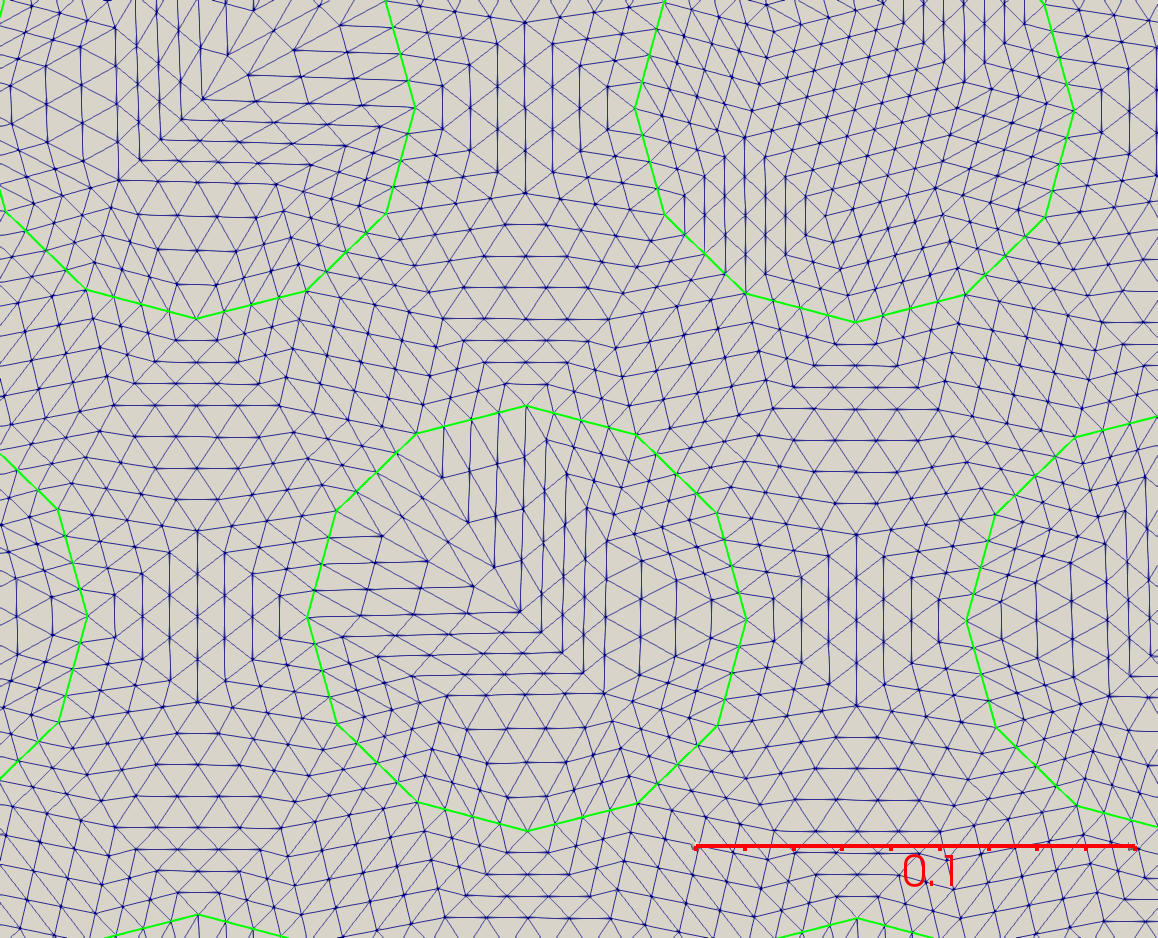}
    \end{subfigure}
    \begin{subfigure}[t]{\textwidth}
        \centering
        \makebox[0.04\textwidth][l]{(b)}
        \includegraphics[width=0.25\textwidth]{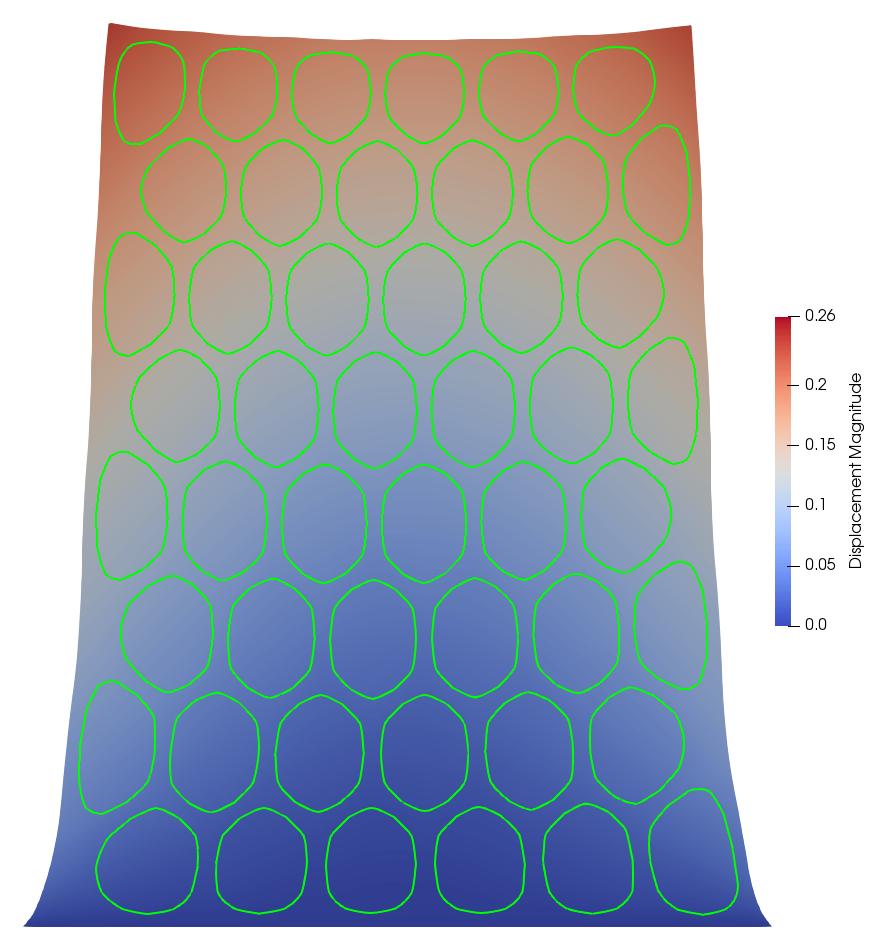}
        %\hfill
        \hspace*{1cm}
        \includegraphics[width=0.25\textwidth]{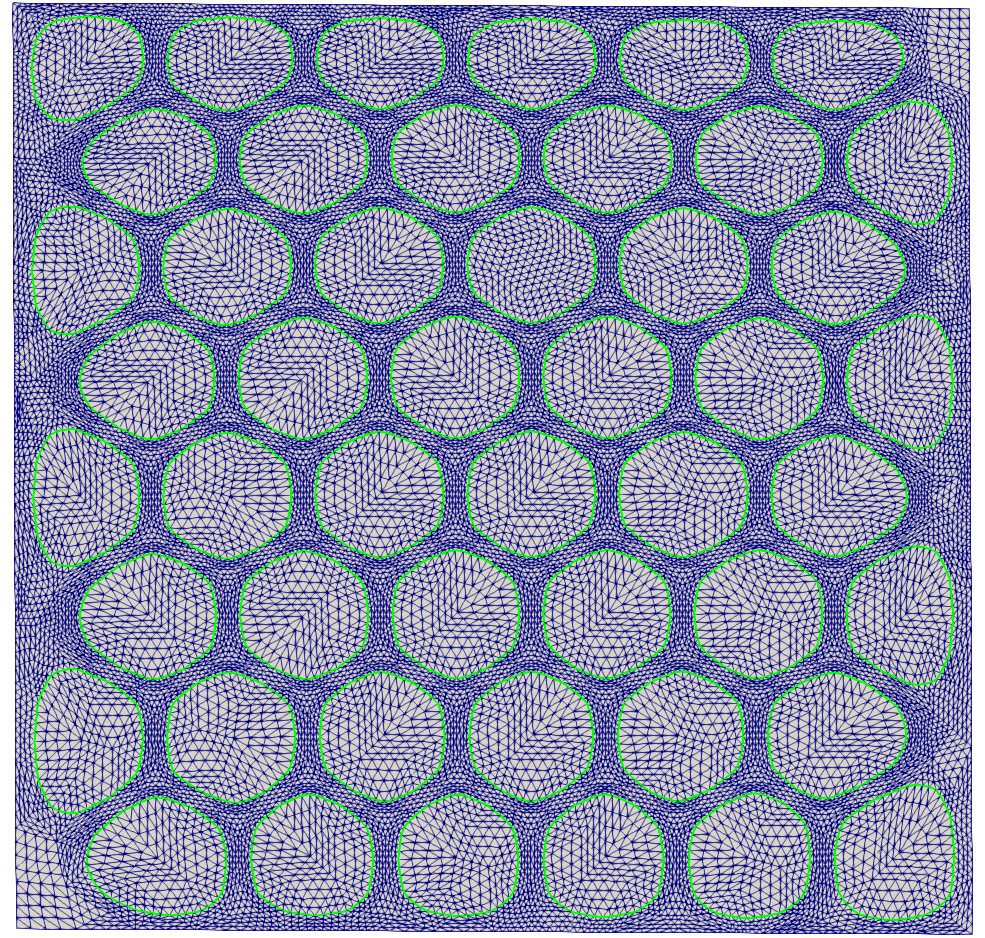}
        \includegraphics[width=0.2\textwidth]{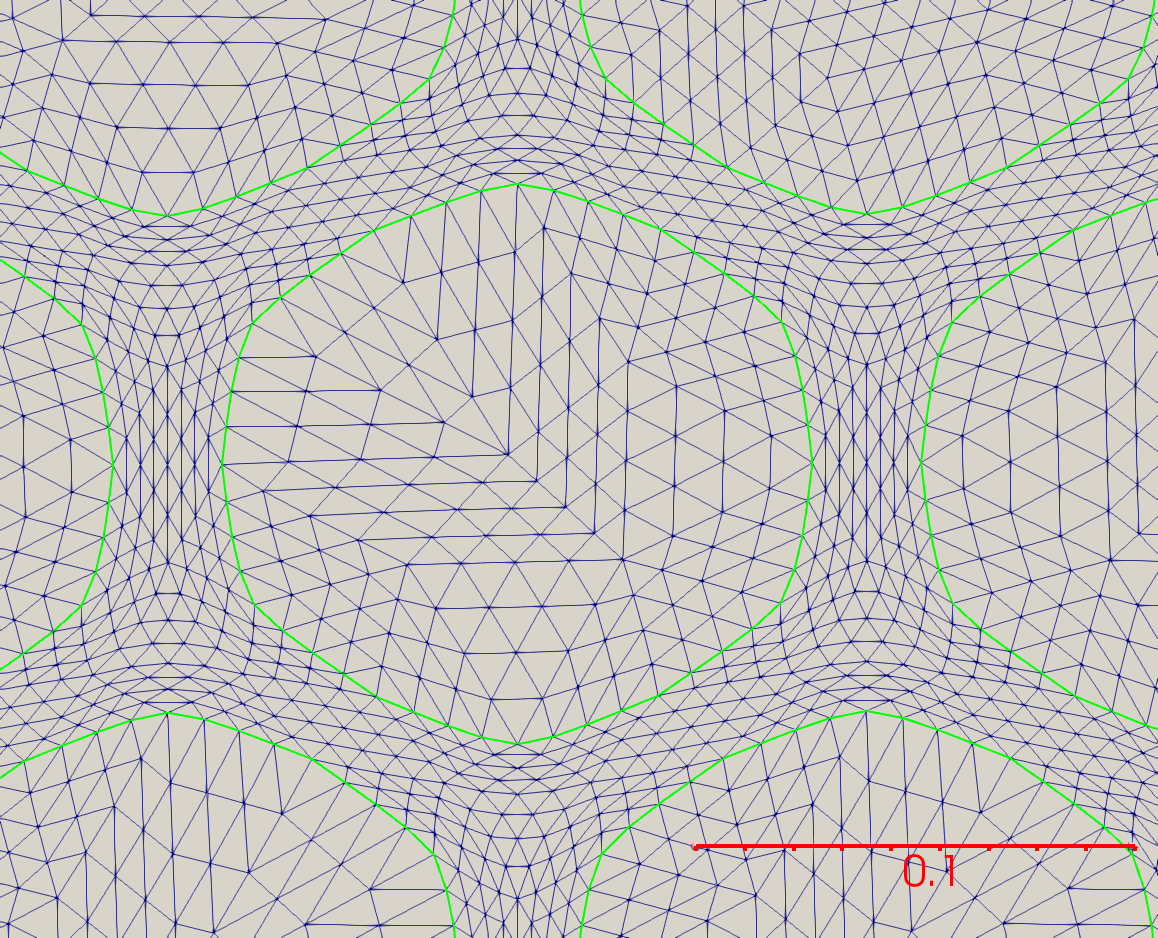}
    \end{subfigure}
    \begin{subfigure}[t]{\textwidth}
        \centering
        \makebox[0.04\textwidth][l]{(c)}
        \includegraphics[width=0.25\textwidth]{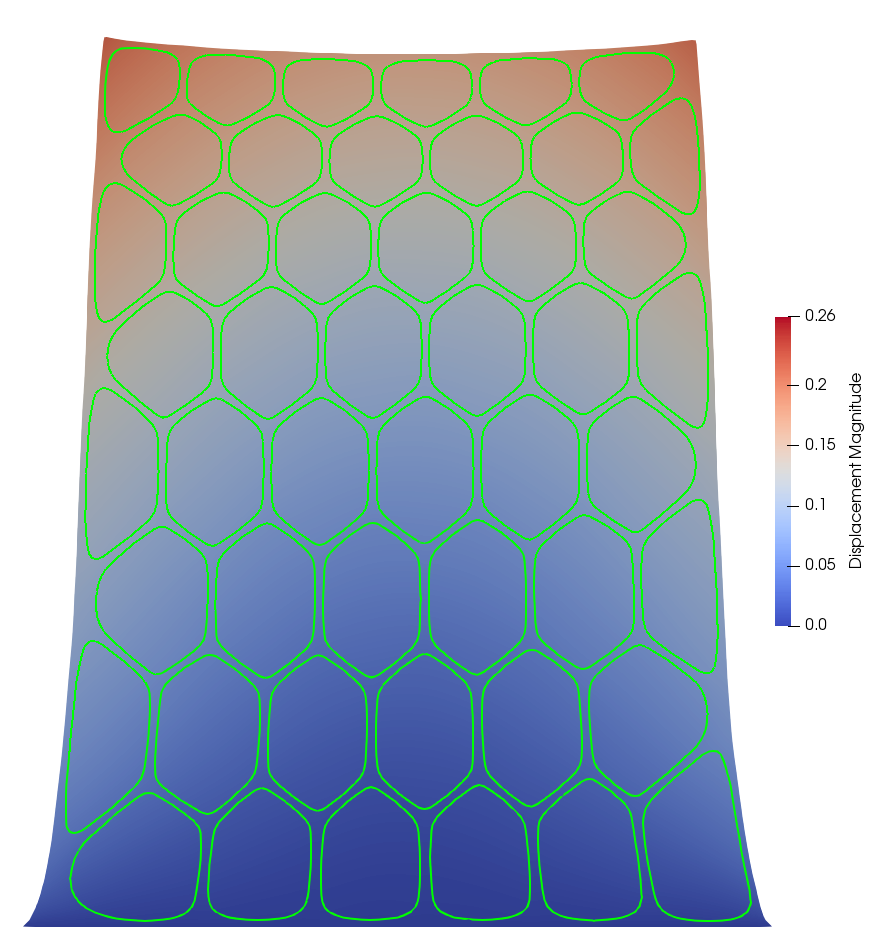}
        %\hfill
        \hspace*{1cm}
        \includegraphics[width=0.25\textwidth]{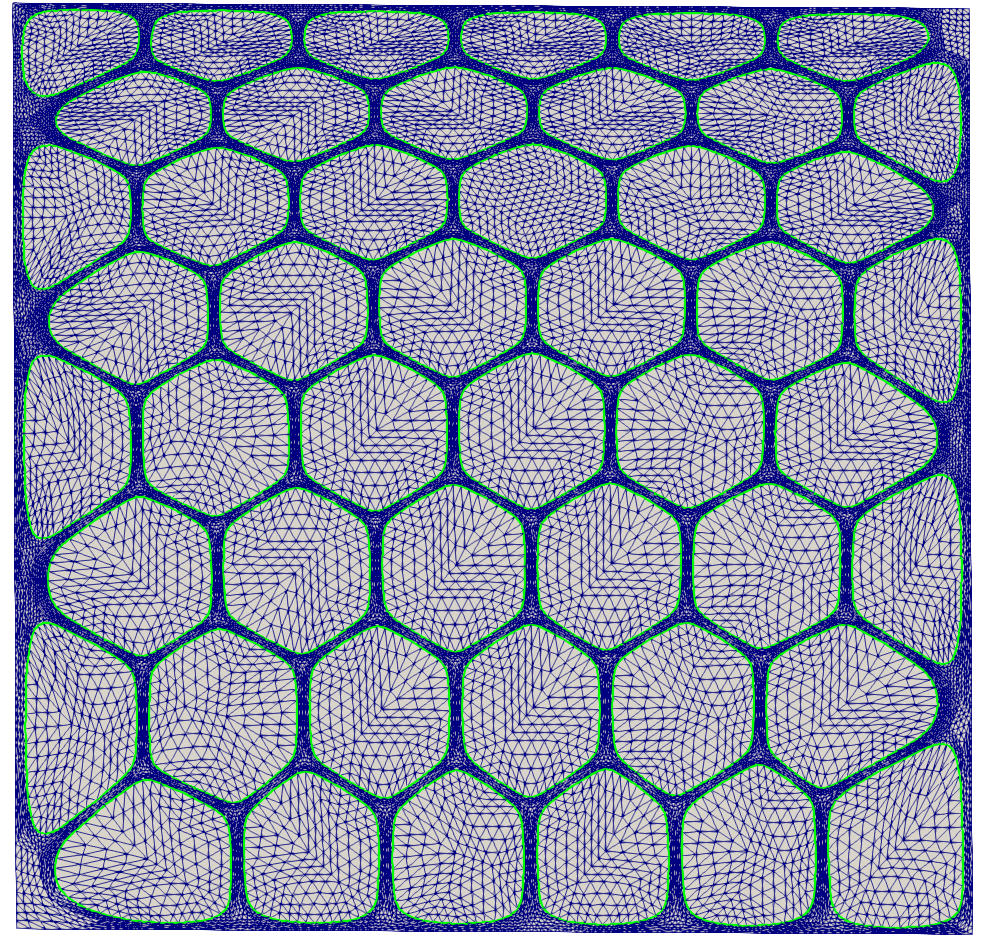}
        \includegraphics[width=0.2\textwidth]{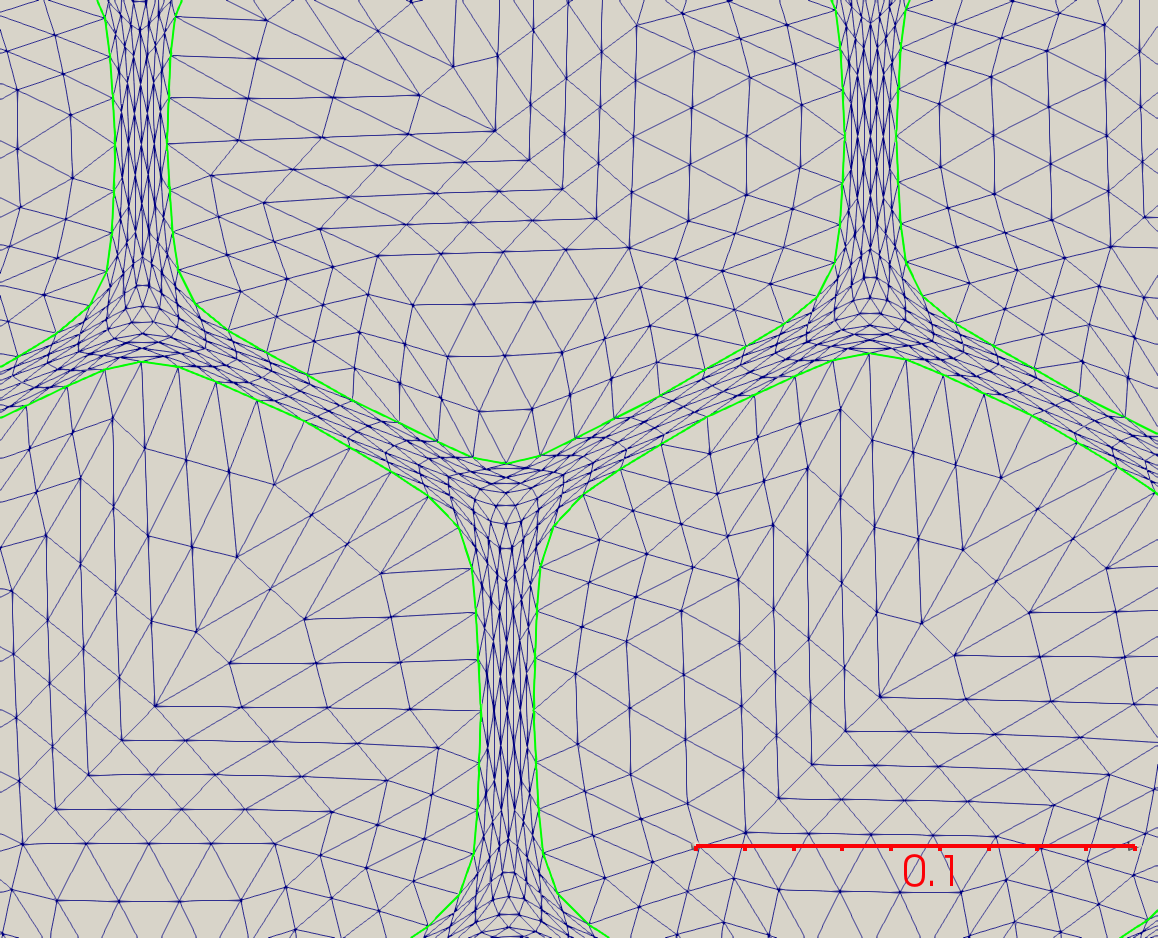}
    \end{subfigure}
    \begin{subfigure}[t]{\textwidth}
        \centering
        \makebox[0.04\textwidth][l]{(d)}
        \includegraphics[width=0.25\textwidth]{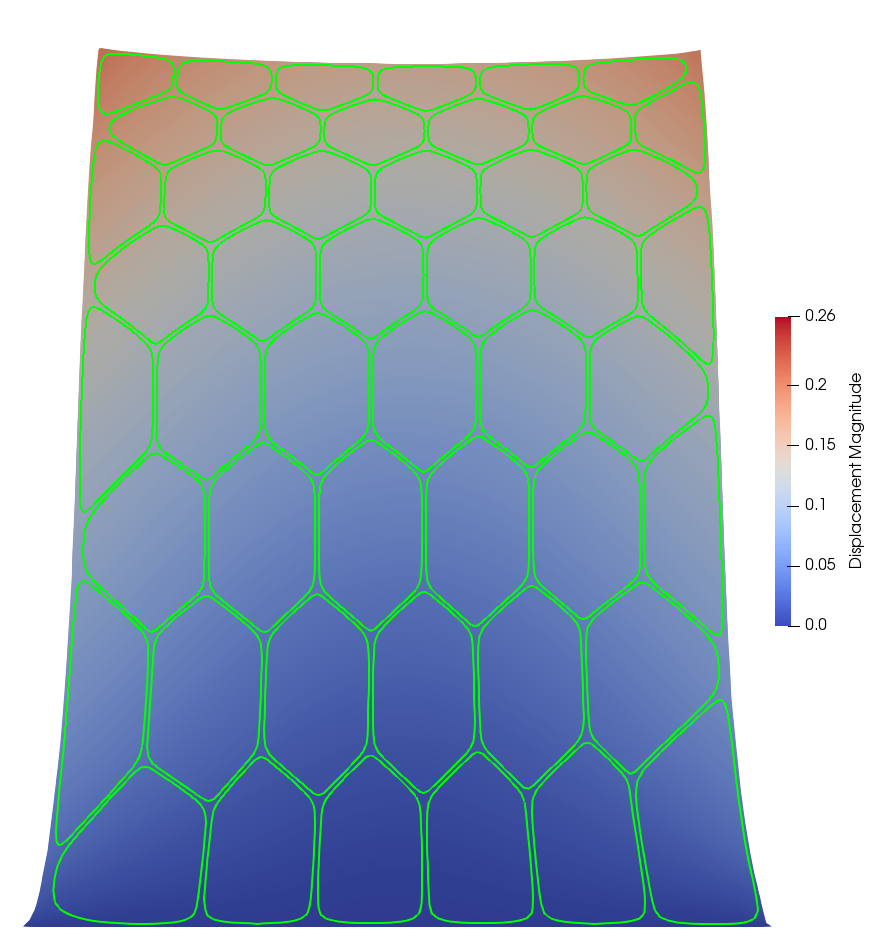}
        %\hfill\\
        \hspace*{1cm}
        \includegraphics[width=0.25\textwidth]{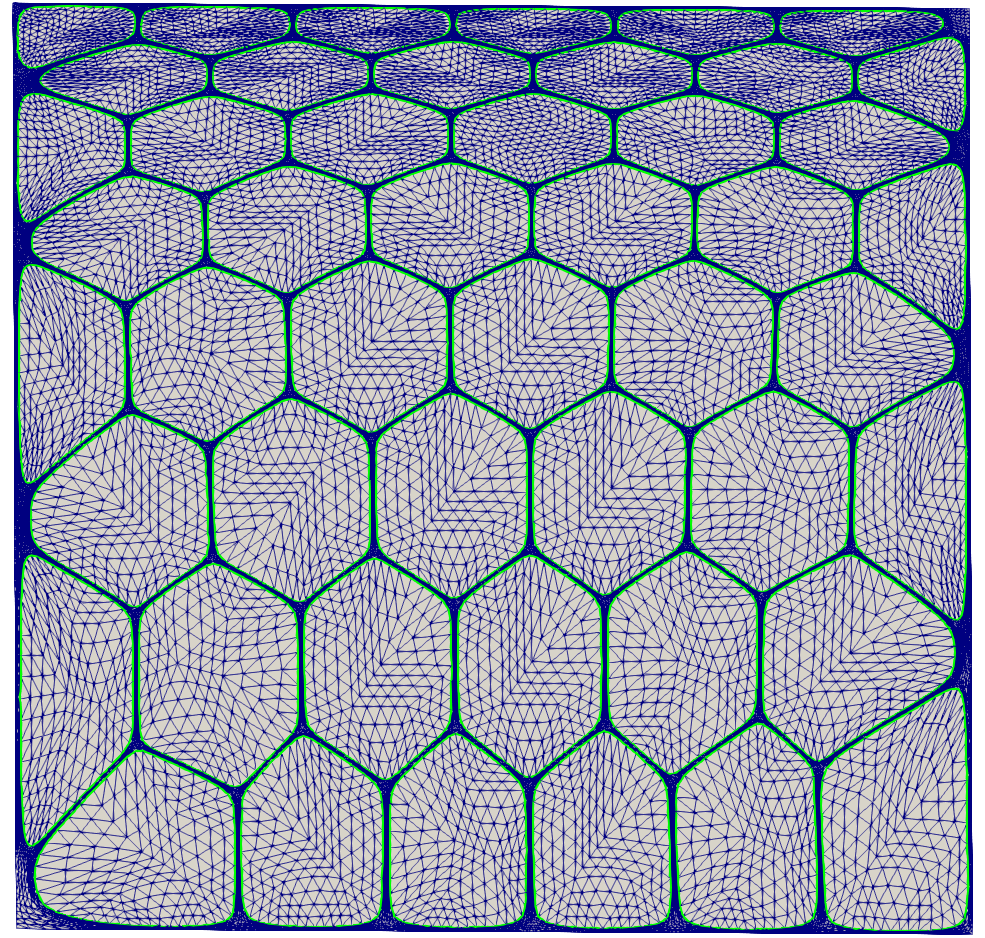}
        \includegraphics[width=0.2\textwidth]{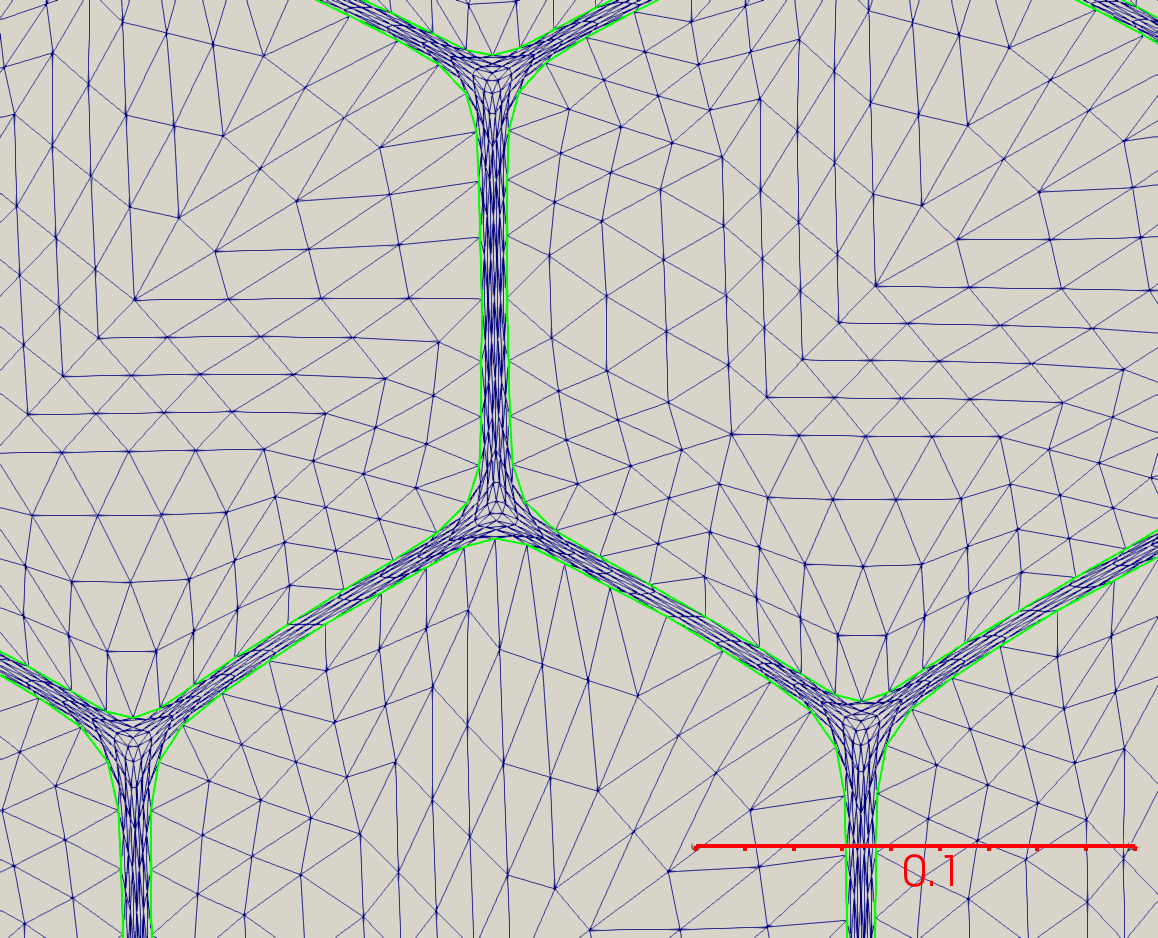}
    \end{subfigure}
    \begin{subfigure}[t]{\textwidth}
        \centering
        \makebox[0.04\textwidth][l]{(e)}
        \includegraphics[width=0.25\textwidth]{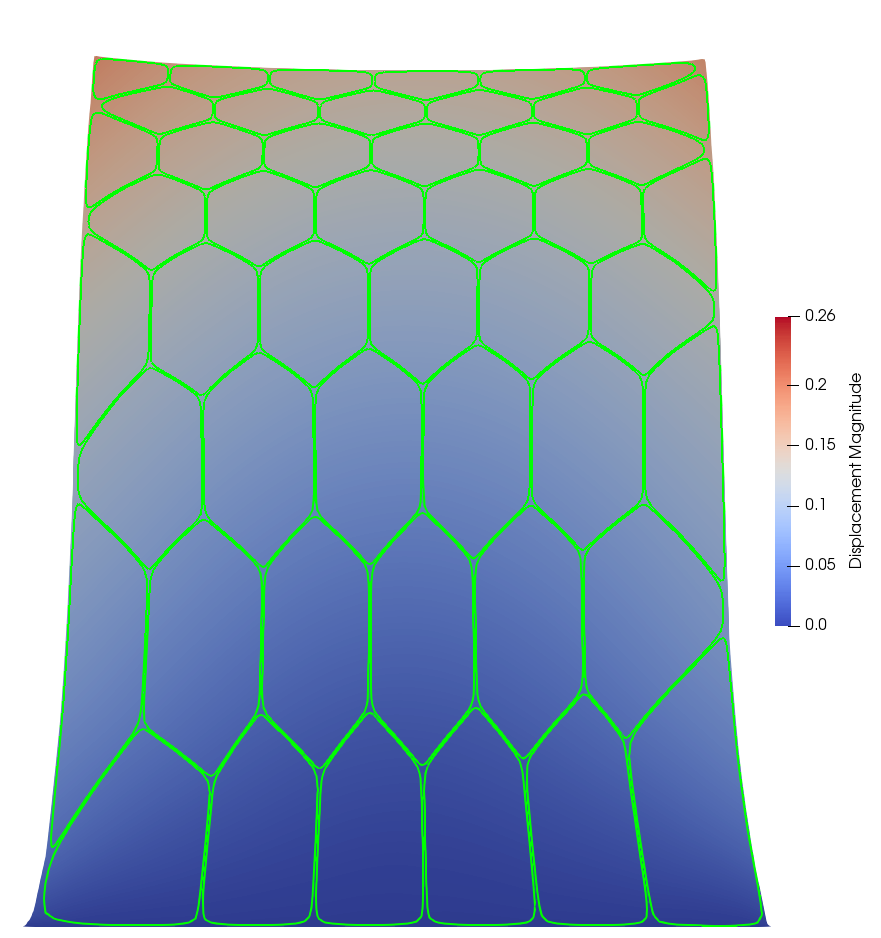}
        %\hfill
        \hspace*{1cm}
        \includegraphics[width=0.25\textwidth]{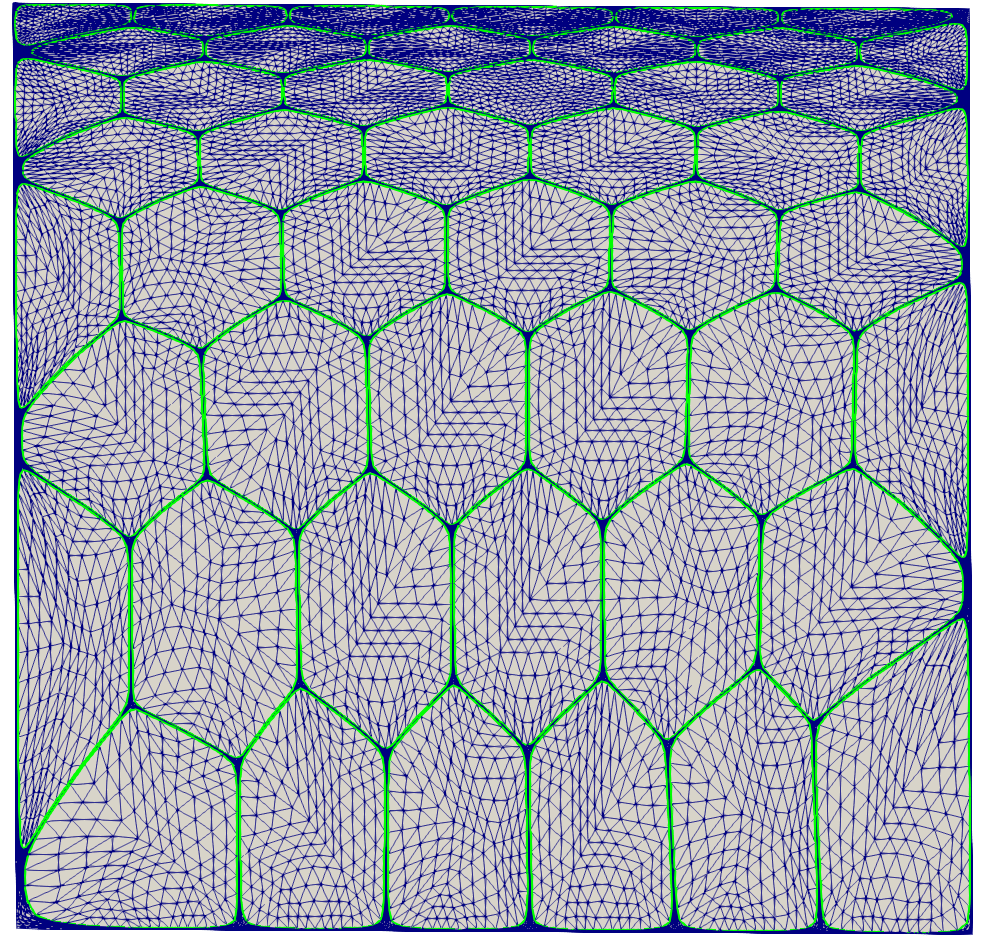}
        \includegraphics[width=0.2\textwidth]{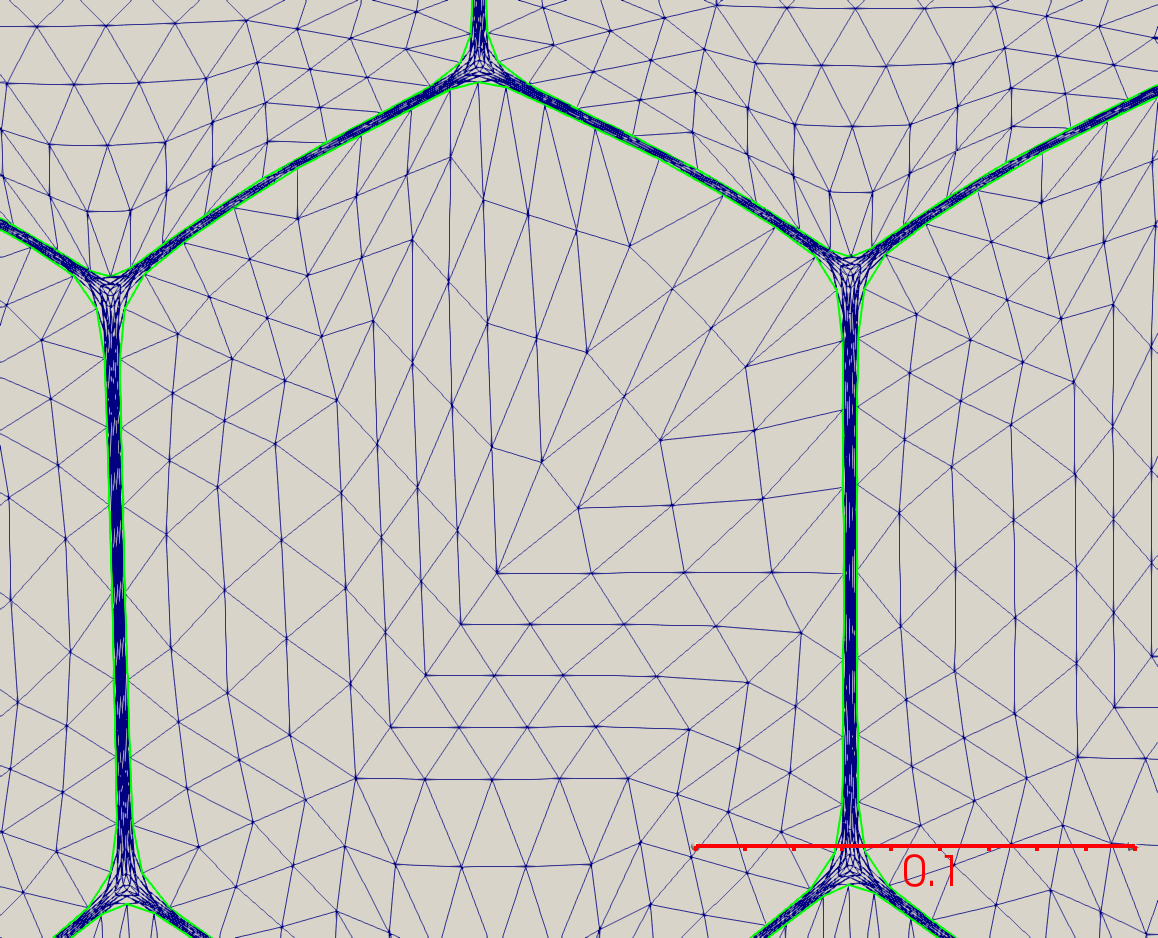}
    \end{subfigure}
    \caption{Shape positions of the cells (green) in the deformed specimen (left) and in the FEM mesh in reference configuration (right) shown for the optimization steps: (a) 1 (b) 25 (c) 50 (d) 75 (e) 100}
    \label{fig:OptimSteps}
\end{figure}

The results of an optimization for $\nuElas=100$, $\nuVol=1$, $\nuPeri=0.01$, $\nuPenalty=5\cdot10^4$, $b=0.001$, $t=1$, and 2 levels of refinement for the FEM mesh is shown in Figure~\ref{fig:OptimSteps}. As solver, we used the Newton scheme without linesearch for the nonlinearity, and the linear system of equations are solved using the BiCGStab method with a geometric multigrid preconditioner employing 3 Block-Jacobi pre- and postsmoothing steps and a LU base solver. 

The simulation results depict the change of the surface shape for selected optimization steps. As seen from the reference configuration, the objective first quickly increases all cells in order to reduce the outer volume, where the penalty term helps to still keep intercell channels open during the optimization. Once the domain is almost completely filled with cells, the optimization turns to an increase of the stiffer cells at the bottom and decreases the volume of the less stiff cells at the top. From the displacement figures, we can observe that the overall displacement diminishes with the optimization which reflects the fact that now the material has become stiffer due to the new size and location of the cell domains.

\begin{figure}[tbp]
    \centering
    \scalebox{0.7}{\input{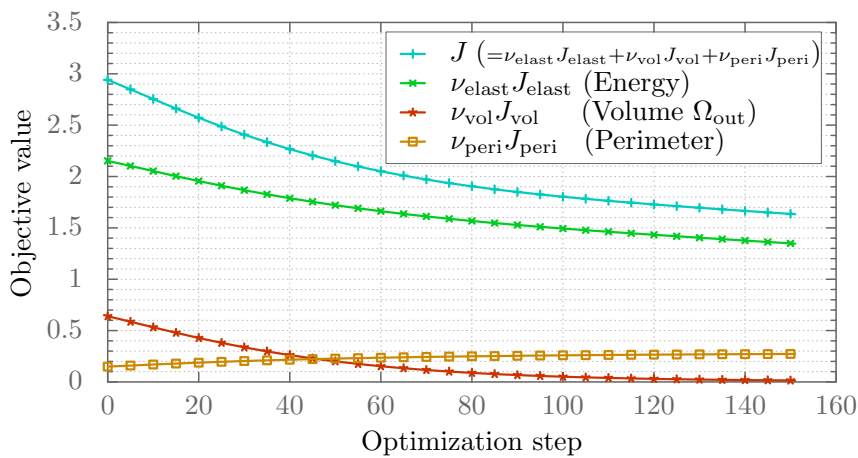}}
    \caption{Value of the objective function for the different components for $b=0.001$}
    \label{fig:ObjectiveValue}
\end{figure}
 
This behavior is also reflected in the behavior of the objective functional and we show the values of the individual components in Figure~\ref{fig:ObjectiveValue}. As can be observed, the computed shape movement is a descent direction for the objective value. Once the outer volume contribution has been removed considerably, the energy contribution is balanced against the surface term by decreasing the size of less stiff cells.

\subsection{Influence of gradient penalization}

\begin{figure}[h!tb]
%    \centering
    \begin{subfigure}[b]{0.24\textwidth}
%        \centering
        \includegraphics[width=0.8\textwidth]{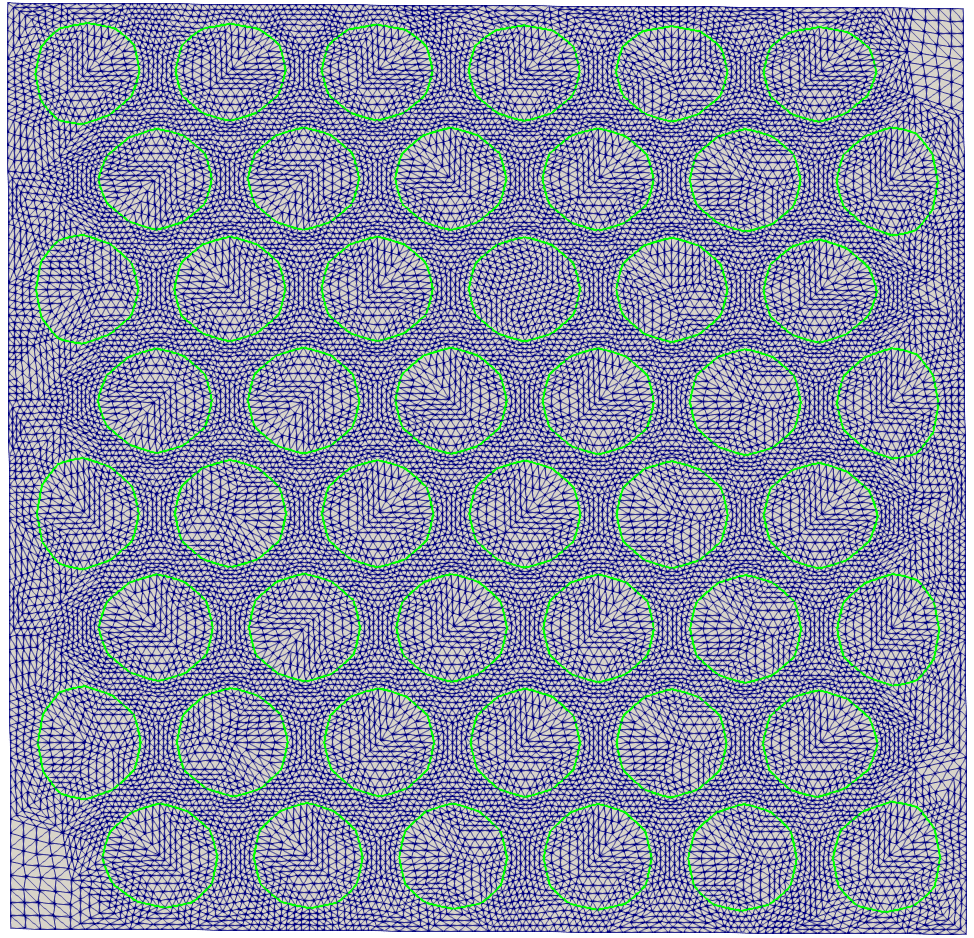}
        \\
        \makebox[0.1\textwidth][l]{(a)}
 %       \vspace{0.2em}
    \end{subfigure}
    \begin{subfigure}[b]{0.24\textwidth}
%        \centering
        \includegraphics[width=0.8\textwidth]{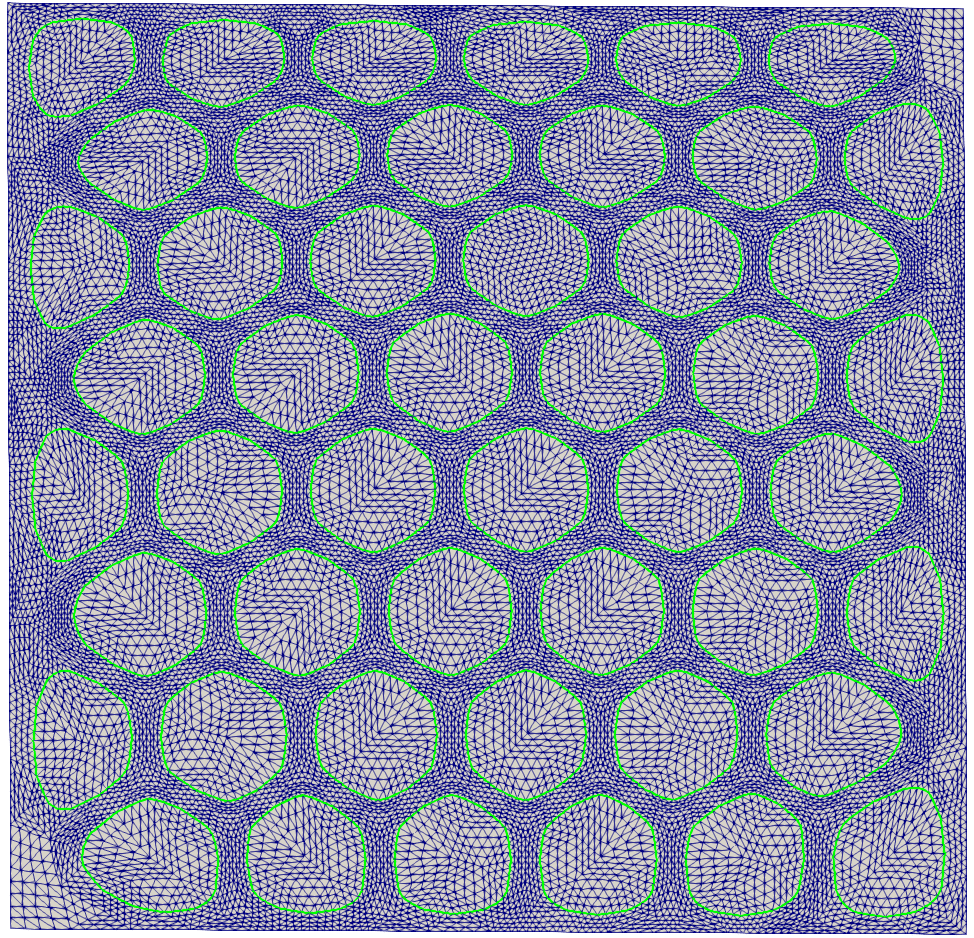}
        \\
        \makebox[0.1\textwidth][l]{(b)}
 %       \vspace{0.2em}
    \end{subfigure}
    \begin{subfigure}[b]{0.24\textwidth}
%        \centering
        \includegraphics[width=0.8\textwidth]{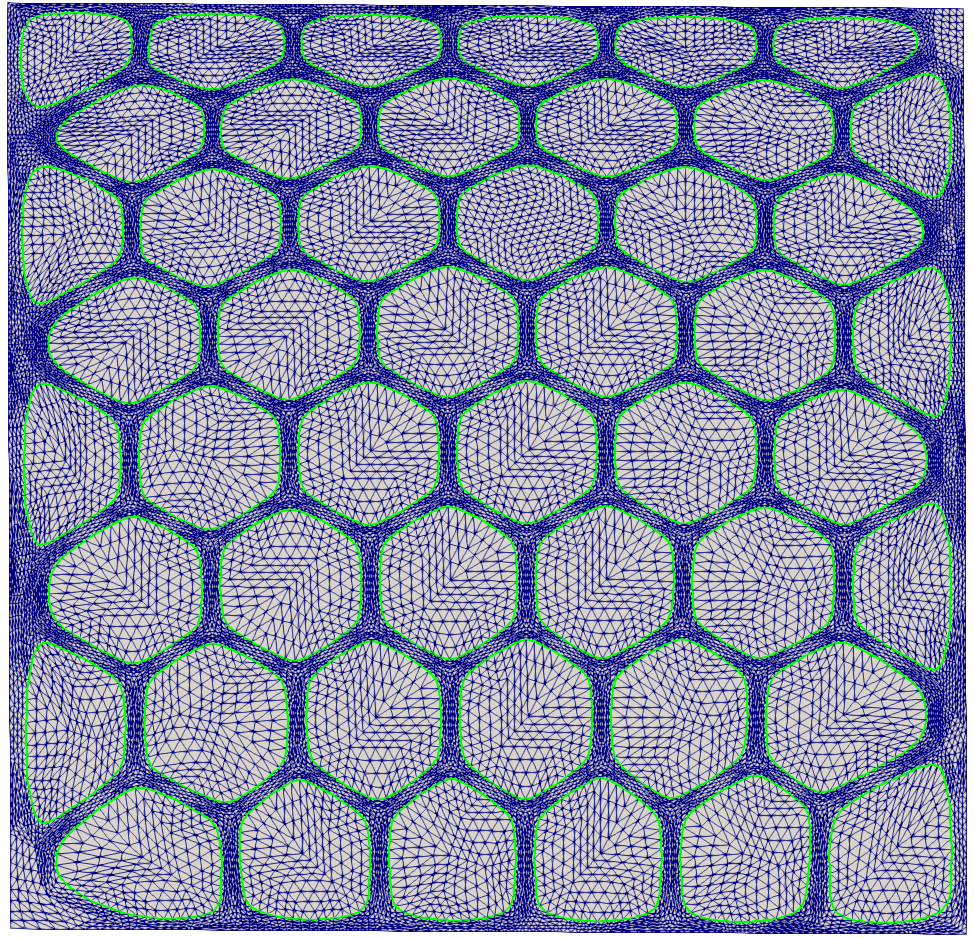}
        \\
        \makebox[0.1\textwidth][l]{(c)}
%        \vspace{0.2em}
    \end{subfigure}
    \begin{subfigure}[b]{0.24\textwidth}
%        \centering
        \makebox[0.1\textwidth][l]{(d)}
        \includegraphics[width=0.8\textwidth]{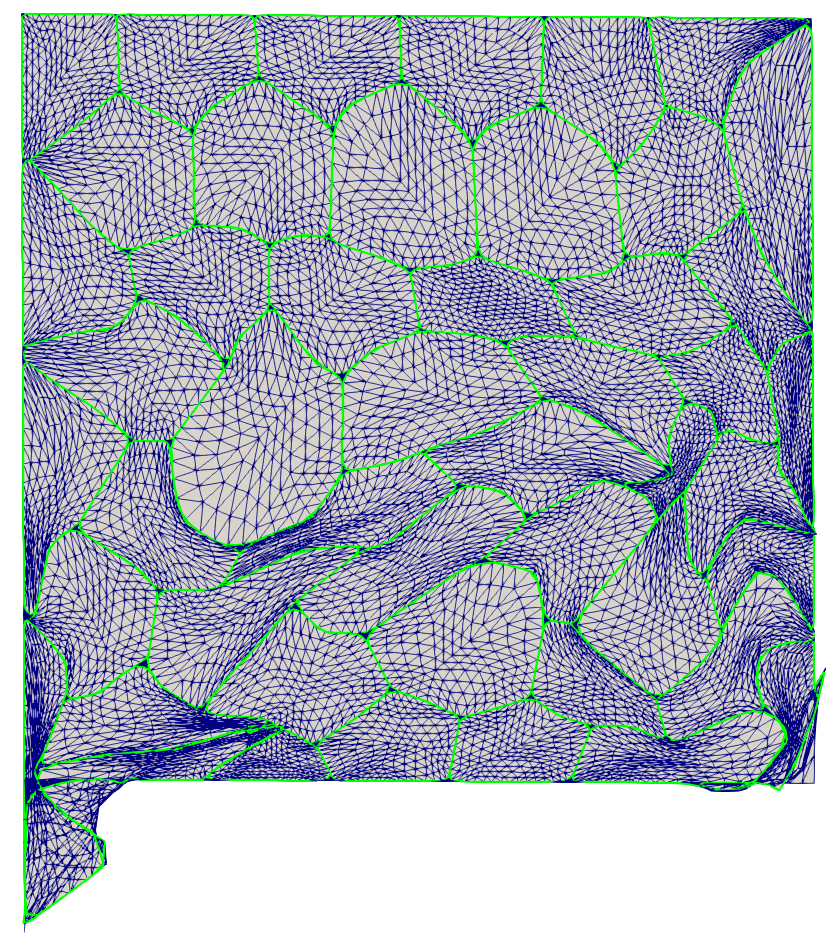}
    \end{subfigure}
    \caption{FEM mesh in reference configuration after 10 optimization steps for different values of $b$: (a) $0.01$, (b) $0.05$, (c) $0.1$; LU enforced solution (d) with mesh penetration for $b=0.0$}
    \label{fig:ComparisionOfB}
\end{figure}

In Figure~\ref{fig:ComparisionOfB}, we show the influence of the penalization of the norm of the gradient with respect to the parameter $b$. With increasing $b$, i.e. allowing larger deformations within each optimization step, the cells are inflated more rapidly. The boundaries of the different shape surfaces therefore approach quicker with two effects: first, the channels between the cells become very close and the mesh is thereby pretty distorted. As a consequence, the optimization stops already after a few (25 for $b=0.1$, 3 for $b=1.0$) steps due to a failure of the iterative solver. Choosing a direct solver like LU factorization in such a case can clearly enforce the computation of a solution, however, it simply results in a penetration of the FEM mesh elements and therefore to meaningless results. Second, due to the oversized (i.e. non-penalized) shape gradient, the optimization does not enter the regime where cell stiffnesses are balanced, but only remains in the inflation part of the optimization.

\subsection{Iteration counts and solver behavior}

For the numerical solution in this example, we have chosen a Newton iteration to resolve the nonlinearity for shape gradient computation which is stopped at an absolute or relative reduction of $10^{-9}$ (whatever applies first). Each linearization within a Newton step is solved by the geometric multigrid preconditioned BiCGstab method down to a relative residuum reduction of $10^{-3}$ which empirically shows to be sufficient to maintain the quadratic reduction rate of the Newton method. The iteration count for the Newton method is restricted to $200$ steps and otherwise the simulation is aborted. For the linear elasticity problem, we apply the geometric multigrid accelerated BiCGStab solver directly and solve down to an absolute or relative reduction of the order $10^{-10}$ while imposing a maximum iteration count of $2000$ steps. 

\begin{figure}[h!tbp]
    \centering
    \begin{subfigure}[t]{0.45\textwidth}
        \centering
        \scalebox{0.49}{\input{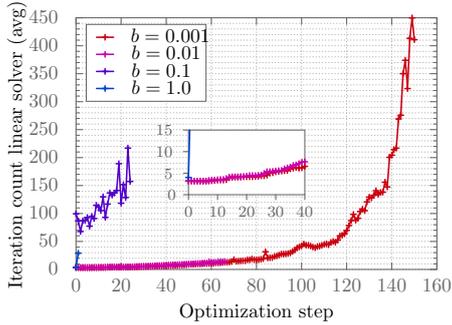}}
        \caption{Average number of linear iteration steps to reduce the residuum by a factor of $10^{-3}$ within a Newton step}
    \end{subfigure}
    \hfill
    \begin{subfigure}[t]{0.45\textwidth}
        \centering
        \scalebox{0.5}{\input{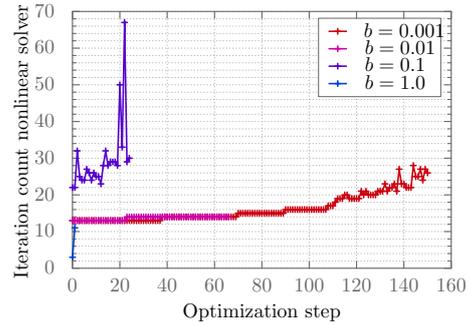}}
        \caption{Number of steps for the nonlinear iteration (Newton method)}
    \end{subfigure}
    \begin{subfigure}[t]{0.45\textwidth}
        \centering
        \scalebox{0.49}{\input{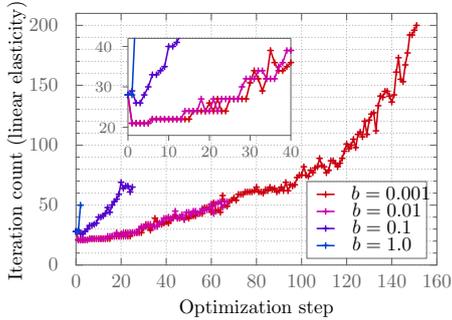}}
       \caption{Number of iterations for the linear solver of the elasticity problem}
    \end{subfigure}
    \hfill
    \begin{subfigure}[t]{0.45\textwidth}
        \centering
        \scalebox{0.5}{\input{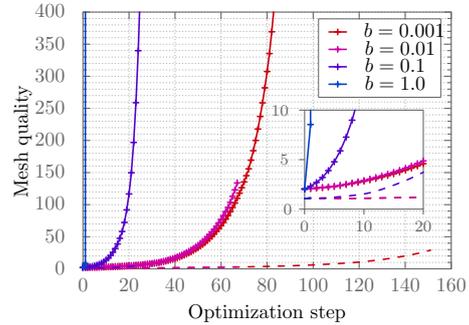}}

        \caption{Mesh quality (ratio between circumsphere and inscribed sphere radius of an element); shown are maximal (solid) and median (dashed) values}
    \end{subfigure}
    \caption{Iteration counts for linear and nonlinear solver (initial mesh 2 times refined)}
    \label{fig:IterationCount}
\end{figure}

The solver behavior for each optimization step is depicted in Figure~\ref{fig:IterationCount}. For the individual optimization steps and different values of $b$, We show the iteration counts for the Newton solver, the average iteration count of the linear solver within a Newton step, and the linear solver iteration count for the elasticity problem. A stopped curve corresponds to a simulation run where the convergence criteria could not be met within the imposed maximal iteration count and the computation was aborted. Generally two observations can be stated: first, with larger values of $b$ (allowing larger deformations for one optimization step), the iteration count for the nonlinear as well as the linear solver quickly increases dramatically and the simulation is stopped. Second, with sufficient penalization, the simulation can run for a reasonable amount of steps. Again, the iteration count increases towards the end of the simulation which we attribute to the distorted finite element mesh. However, with penalization and up to 100 optimization steps (see Fig.~\ref{fig:OptimSteps} for a visualization), the iteration count for the Newton method roughly remains constant (approx. $15$ steps for each optimization step), and the average linear iteration count for the relative reduction within a newton step increases from $3$ to a moderate number of $40$ steps. For the linear elasticity solver, the iteration count increases from $20$ for the nicely shaped mesh at the beginning to roughly $75$ steps for the distorted mesh at optimization step $100$.

\subsection{Parallelization aspects}

For the proposed method, the topology of the finite element mesh is fixed and not altered during the simulation. This clearly has the drawback that element shapes become unpleasant, and with the deteriorating mesh quality also a reduction in FE solution quality as well as in the solver efficiency has to be expected. This has been observed for the iteration counts in our examples, however, the introduced penalization helps to mitigate this behavior and to allow for a reasonable optimization. 
On the other hand, the fixed FE mesh topology has an important advantage: with respect to parallelization, the mesh has not to be rebalanced or redistributed among involved processes. An initial, well-balanced partition of the finite element mesh onto the processes will remain for the whole optimization range and this reasoning applies also to the load balance of coarse mesh hierarchies. Therefore, using available and established linear solver implementations, a good scaling is to be expected once a good initial load balancing can be achieved. 

For this study, we use the parallel geometric multigrid approach with hierarchically distributed meshes \cite{Reiter2014gmg}. This approach has shown close to optimal scalability on hundred thousands of cores. Our mesh partitioning is based on the ParMETIS \cite{ParMetis} software. Here, we perform a strong and weak scaling study for the optimization components and measure accumulated times for the first 10 optimization steps. We use the BiCGStab method preconditioned with the geometric multigrid using a Block-Jacobi smoother with 3 pre- and postsmoothing steps, V-cycle, assembled coarse grid matrices, and a gathered LU base solver. All measurements have been performed on the \emph{Statik\&Dynamik}-Cluster at the Department of Civil and Environmental Engineering of the Ruhr University Bochum. The Cluster consists of 92 nodes, each with two Intel Xeon Skylake Gold 6148 (40 cores per node) and 192GB memory, i.e. a maximum of 3680 cores. Our study tries to fully utilize the allocated nodes and thus starts at 40 cores and increases the core count by a factor of 4 (refining the 2d mesh with each increase for the weak scaling study) up to 2560 cores. In order to account for jitter, we have performed five runs for each measurement and present the arithmetic mean. 

The results are presented in Figure~\ref{fig:StrongScaling} for the strong scaling, and in Figure~\ref{fig:WeakScaling} for the weak scaling study. We present wallclock times (accumulated over 10 optimization steps) and the gained speedup relative to 40 cores. Both studies show a good scalability and speedup, and the number of Newton steps remains almost constant even for the weak scaling study. A surprising fact is the increase in the iteration count of the linear solver. While this is to be expected for the weak scaling (a slightly different numerical problem is solved with each mesh refinement), this should not be the case for the strong scaling study. Note that we employ a Jacobi smoother and the mathematical algorithm is thereby exactly the same even for parallel computations. A detailed investigation revealed that small differences in round-off errors for different parallel partitions accumulate for these simulations and explain the slightly different iteration counts.  

\begin{figure}[h!tbp]

  \centering
  \begin{subfigure}[t]{\textwidth}
    \centering
    \scalebox{0.48}{\begin{tikzpicture}[gnuplot]
%% generated with GNUPLOT 5.2p6 (Lua 5.3; terminal rev. Nov 2018, script rev. 107)
%% Fr  8 Mär 17:47:29 2019
\tikzset{every node/.append style={scale=1.50}}
\path (0.000,0.000) rectangle (12.500,8.750);
\gpcolor{rgb color={0.702,0.702,0.702}}
\gpsetlinetype{gp lt axes}
\gpsetdashtype{gp dt axes}
\gpsetlinewidth{1.00}
\draw[gp path] (2.115,1.852)--(11.671,1.852);
\gpcolor{rgb color={0.400,0.400,0.400}}
\gpsetlinetype{gp lt border}
\gpsetdashtype{gp dt solid}
\gpsetlinewidth{2.00}
\draw[gp path] (2.115,1.852)--(2.048,1.852);
\gpcolor{rgb color={0.702,0.702,0.702}}
\gpsetlinetype{gp lt axes}
\gpsetdashtype{gp dt axes}
\gpsetlinewidth{1.00}
\draw[gp path] (2.115,2.091)--(11.671,2.091);
\gpcolor{rgb color={0.400,0.400,0.400}}
\gpsetlinetype{gp lt border}
\gpsetdashtype{gp dt solid}
\gpsetlinewidth{2.00}
\draw[gp path] (2.115,2.091)--(2.048,2.091);
\gpcolor{rgb color={0.702,0.702,0.702}}
\gpsetlinetype{gp lt axes}
\gpsetdashtype{gp dt axes}
\gpsetlinewidth{1.00}
\draw[gp path] (2.115,2.283)--(11.671,2.283);
\gpcolor{rgb color={0.400,0.400,0.400}}
\gpsetlinetype{gp lt border}
\gpsetdashtype{gp dt solid}
\gpsetlinewidth{2.00}
\draw[gp path] (2.115,2.283)--(1.980,2.283);
\node[gp node right,font={\fontsize{8.0pt}{9.6pt}\selectfont}] at (1.980,2.283) {1};
\gpcolor{rgb color={0.702,0.702,0.702}}
\gpsetlinetype{gp lt axes}
\gpsetdashtype{gp dt axes}
\gpsetlinewidth{1.00}
\draw[gp path] (2.115,3.277)--(11.671,3.277);
\gpcolor{rgb color={0.400,0.400,0.400}}
\gpsetlinetype{gp lt border}
\gpsetdashtype{gp dt solid}
\gpsetlinewidth{2.00}
\draw[gp path] (2.115,3.277)--(2.048,3.277);
\gpcolor{rgb color={0.702,0.702,0.702}}
\gpsetlinetype{gp lt axes}
\gpsetdashtype{gp dt axes}
\gpsetlinewidth{1.00}
\draw[gp path] (2.115,3.757)--(11.671,3.757);
\gpcolor{rgb color={0.400,0.400,0.400}}
\gpsetlinetype{gp lt border}
\gpsetdashtype{gp dt solid}
\gpsetlinewidth{2.00}
\draw[gp path] (2.115,3.757)--(2.048,3.757);
\gpcolor{rgb color={0.702,0.702,0.702}}
\gpsetlinetype{gp lt axes}
\gpsetdashtype{gp dt axes}
\gpsetlinewidth{1.00}
\draw[gp path] (2.115,4.076)--(11.671,4.076);
\gpcolor{rgb color={0.400,0.400,0.400}}
\gpsetlinetype{gp lt border}
\gpsetdashtype{gp dt solid}
\gpsetlinewidth{2.00}
\draw[gp path] (2.115,4.076)--(2.048,4.076);
\gpcolor{rgb color={0.702,0.702,0.702}}
\gpsetlinetype{gp lt axes}
\gpsetdashtype{gp dt axes}
\gpsetlinewidth{1.00}
\draw[gp path] (2.115,4.316)--(11.671,4.316);
\gpcolor{rgb color={0.400,0.400,0.400}}
\gpsetlinetype{gp lt border}
\gpsetdashtype{gp dt solid}
\gpsetlinewidth{2.00}
\draw[gp path] (2.115,4.316)--(2.048,4.316);
\gpcolor{rgb color={0.702,0.702,0.702}}
\gpsetlinetype{gp lt axes}
\gpsetdashtype{gp dt axes}
\gpsetlinewidth{1.00}
\draw[gp path] (2.115,4.507)--(11.671,4.507);
\gpcolor{rgb color={0.400,0.400,0.400}}
\gpsetlinetype{gp lt border}
\gpsetdashtype{gp dt solid}
\gpsetlinewidth{2.00}
\draw[gp path] (2.115,4.507)--(1.980,4.507);
\node[gp node right,font={\fontsize{8.0pt}{9.6pt}\selectfont}] at (1.888,4.507) {10};
\gpcolor{rgb color={0.702,0.702,0.702}}
\gpsetlinetype{gp lt axes}
\gpsetdashtype{gp dt axes}
\gpsetlinewidth{1.00}
\draw[gp path] (2.115,5.502)--(11.671,5.502);
\gpcolor{rgb color={0.400,0.400,0.400}}
\gpsetlinetype{gp lt border}
\gpsetdashtype{gp dt solid}
\gpsetlinewidth{2.00}
\draw[gp path] (2.115,5.502)--(2.048,5.502);
\gpcolor{rgb color={0.702,0.702,0.702}}
\gpsetlinetype{gp lt axes}
\gpsetdashtype{gp dt axes}
\gpsetlinewidth{1.00}
\draw[gp path] (2.115,5.982)--(11.671,5.982);
\gpcolor{rgb color={0.400,0.400,0.400}}
\gpsetlinetype{gp lt border}
\gpsetdashtype{gp dt solid}
\gpsetlinewidth{2.00}
\draw[gp path] (2.115,5.982)--(2.048,5.982);
\gpcolor{rgb color={0.702,0.702,0.702}}
\gpsetlinetype{gp lt axes}
\gpsetdashtype{gp dt axes}
\gpsetlinewidth{1.00}
\draw[gp path] (2.115,6.301)--(11.671,6.301);
\gpcolor{rgb color={0.400,0.400,0.400}}
\gpsetlinetype{gp lt border}
\gpsetdashtype{gp dt solid}
\gpsetlinewidth{2.00}
\draw[gp path] (2.115,6.301)--(2.048,6.301);
\gpcolor{rgb color={0.702,0.702,0.702}}
\gpsetlinetype{gp lt axes}
\gpsetdashtype{gp dt axes}
\gpsetlinewidth{1.00}
\draw[gp path] (2.115,6.540)--(8.179,6.540);
\draw[gp path] (11.395,6.540)--(11.671,6.540);
\gpcolor{rgb color={0.400,0.400,0.400}}
\gpsetlinetype{gp lt border}
\gpsetdashtype{gp dt solid}
\gpsetlinewidth{2.00}
\draw[gp path] (2.115,6.540)--(2.048,6.540);
\gpcolor{rgb color={0.702,0.702,0.702}}
\gpsetlinetype{gp lt axes}
\gpsetdashtype{gp dt axes}
\gpsetlinewidth{1.00}
\draw[gp path] (2.115,6.732)--(8.179,6.732);
\draw[gp path] (11.395,6.732)--(11.671,6.732);
\gpcolor{rgb color={0.400,0.400,0.400}}
\gpsetlinetype{gp lt border}
\gpsetdashtype{gp dt solid}
\gpsetlinewidth{2.00}
\draw[gp path] (2.115,6.732)--(1.980,6.732);
\node[gp node right,font={\fontsize{8.0pt}{9.6pt}\selectfont}] at (1.888,6.732) {100};
\gpcolor{rgb color={0.702,0.702,0.702}}
\gpsetlinetype{gp lt axes}
\gpsetdashtype{gp dt axes}
\gpsetlinewidth{1.00}
\draw[gp path] (2.115,7.727)--(8.179,7.727);
\draw[gp path] (11.395,7.727)--(11.671,7.727);
\gpcolor{rgb color={0.400,0.400,0.400}}
\gpsetlinetype{gp lt border}
\gpsetdashtype{gp dt solid}
\gpsetlinewidth{2.00}
\draw[gp path] (2.115,7.727)--(2.048,7.727);
\gpcolor{rgb color={0.702,0.702,0.702}}
\gpsetlinetype{gp lt axes}
\gpsetdashtype{gp dt axes}
\gpsetlinewidth{1.00}
\draw[gp path] (2.115,8.206)--(11.671,8.206);
\gpcolor{rgb color={0.400,0.400,0.400}}
\gpsetlinetype{gp lt border}
\gpsetdashtype{gp dt solid}
\gpsetlinewidth{2.00}
\draw[gp path] (2.115,8.206)--(2.048,8.206);
\gpcolor{rgb color={0.702,0.702,0.702}}
\gpsetlinetype{gp lt axes}
\gpsetdashtype{gp dt axes}
\gpsetlinewidth{1.00}
\draw[gp path] (2.682,1.613)--(2.682,8.287);
\gpcolor{rgb color={0.400,0.400,0.400}}
\gpsetlinetype{gp lt border}
\gpsetdashtype{gp dt solid}
\gpsetlinewidth{2.00}
\draw[gp path] (2.682,1.613)--(2.682,1.478);
\node[gp node center,font={\fontsize{8.0pt}{9.6pt}\selectfont}] at (2.682,1.016) {40};
\gpcolor{rgb color={0.702,0.702,0.702}}
\gpsetlinetype{gp lt axes}
\gpsetdashtype{gp dt axes}
\gpsetlinewidth{1.00}
\draw[gp path] (5.412,1.613)--(5.412,8.287);
\gpcolor{rgb color={0.400,0.400,0.400}}
\gpsetlinetype{gp lt border}
\gpsetdashtype{gp dt solid}
\gpsetlinewidth{2.00}
\draw[gp path] (5.412,1.613)--(5.412,1.478);
\node[gp node center,font={\fontsize{8.0pt}{9.6pt}\selectfont}] at (5.412,1.016) {160};
\gpcolor{rgb color={0.702,0.702,0.702}}
\gpsetlinetype{gp lt axes}
\gpsetdashtype{gp dt axes}
\gpsetlinewidth{1.00}
\draw[gp path] (8.142,1.613)--(8.142,8.287);
\gpcolor{rgb color={0.400,0.400,0.400}}
\gpsetlinetype{gp lt border}
\gpsetdashtype{gp dt solid}
\gpsetlinewidth{2.00}
\draw[gp path] (8.142,1.613)--(8.142,1.478);
\node[gp node center,font={\fontsize{8.0pt}{9.6pt}\selectfont}] at (8.142,1.016) {640};
\gpcolor{rgb color={0.702,0.702,0.702}}
\gpsetlinetype{gp lt axes}
\gpsetdashtype{gp dt axes}
\gpsetlinewidth{1.00}
\draw[gp path] (10.872,1.613)--(10.872,6.490);
\draw[gp path] (10.872,8.107)--(10.872,8.287);
\gpcolor{rgb color={0.400,0.400,0.400}}
\gpsetlinetype{gp lt border}
\gpsetdashtype{gp dt solid}
\gpsetlinewidth{2.00}
\draw[gp path] (10.872,1.613)--(10.872,1.478);
\node[gp node center,font={\fontsize{8.0pt}{9.6pt}\selectfont}] at (10.872,1.016) {2560};
\draw[gp path] (2.115,8.287)--(2.115,1.613)--(11.671,1.613)--(11.671,8.287)--cycle;
\gpcolor{color=gp lt color border}
\node[gp node center,rotate=-270] at (0.874,4.950) {Time [s]};
\node[gp node center] at (6.893,0.400) {Number of processes};
\gpcolor{rgb color={0.400,0.400,0.400}}
\gpsetlinewidth{1.00}
\draw[gp path] (8.179,6.490)--(8.179,8.107)--(11.395,8.107)--(11.395,6.490)--cycle;
\gpcolor{color=gp lt color border}
\node[gp node left] at (9.187,7.760) {$T_{\text{ass}}$};
\gpcolor{rgb color={0.800,0.000,0.063}}
\gpsetlinewidth{3.00}
\draw[gp path] (8.455,7.760)--(8.911,7.760);
\draw[gp path] (2.682,6.423)--(5.412,5.224)--(8.142,4.002)--(10.872,2.734);
\gpsetpointsize{5.20}
\gppoint{gp mark 1}{(2.682,6.423)}
\gppoint{gp mark 1}{(5.412,5.224)}
\gppoint{gp mark 1}{(8.142,4.002)}
\gppoint{gp mark 1}{(10.872,2.734)}
\gppoint{gp mark 1}{(8.683,7.760)}
\gpcolor{color=gp lt color border}
\node[gp node left] at (9.187,7.298) {$T_{\text{init}}$};
\gpcolor{rgb color={0.800,0.000,0.663}}
\draw[gp path] (8.455,7.298)--(8.911,7.298);
\draw[gp path] (2.682,5.495)--(5.412,4.284)--(8.142,3.103)--(10.872,2.104);
\gppoint{gp mark 2}{(2.682,5.495)}
\gppoint{gp mark 2}{(5.412,4.284)}
\gppoint{gp mark 2}{(8.142,3.103)}
\gppoint{gp mark 2}{(10.872,2.104)}
\gppoint{gp mark 2}{(8.683,7.298)}
\gpcolor{color=gp lt color border}
\node[gp node left] at (9.187,6.836) {$T_{\text{solve}}$};
\gpcolor{rgb color={0.329,0.000,0.800}}
\draw[gp path] (8.455,6.836)--(8.911,6.836);
\draw[gp path] (2.682,7.383)--(5.412,6.015)--(8.142,4.687)--(10.872,3.613);
\gppoint{gp mark 3}{(2.682,7.383)}
\gppoint{gp mark 3}{(5.412,6.015)}
\gppoint{gp mark 3}{(8.142,4.687)}
\gppoint{gp mark 3}{(10.872,3.613)}
\gppoint{gp mark 3}{(8.683,6.836)}
%% coordinates of the plot area
\gpdefrectangularnode{gp plot 1}{\pgfpoint{2.115cm}{1.613cm}}{\pgfpoint{11.671cm}{8.287cm}}
\end{tikzpicture}
%% gnuplot variables
}
    \scalebox{0.48}{\begin{tikzpicture}[gnuplot]
%% generated with GNUPLOT 5.2p6 (Lua 5.3; terminal rev. Nov 2018, script rev. 107)
%% Fr  8 Mär 17:47:29 2019
\tikzset{every node/.append style={scale=1.50}}
\path (0.000,0.000) rectangle (12.500,8.750);
\gpcolor{rgb color={0.702,0.702,0.702}}
\gpsetlinetype{gp lt axes}
\gpsetdashtype{gp dt axes}
\gpsetlinewidth{1.00}
\draw[gp path] (1.839,2.525)--(11.671,2.525);
\gpcolor{rgb color={0.400,0.400,0.400}}
\gpsetlinetype{gp lt border}
\gpsetdashtype{gp dt solid}
\gpsetlinewidth{2.00}
\draw[gp path] (1.839,2.525)--(1.704,2.525);
\node[gp node right,font={\fontsize{8.0pt}{9.6pt}\selectfont}] at (1.704,2.525) {1};
\gpcolor{rgb color={0.702,0.702,0.702}}
\gpsetlinetype{gp lt axes}
\gpsetdashtype{gp dt axes}
\gpsetlinewidth{1.00}
\draw[gp path] (1.839,3.676)--(11.671,3.676);
\gpcolor{rgb color={0.400,0.400,0.400}}
\gpsetlinetype{gp lt border}
\gpsetdashtype{gp dt solid}
\gpsetlinewidth{2.00}
\draw[gp path] (1.839,3.676)--(1.772,3.676);
\gpcolor{rgb color={0.702,0.702,0.702}}
\gpsetlinetype{gp lt axes}
\gpsetdashtype{gp dt axes}
\gpsetlinewidth{1.00}
\draw[gp path] (1.839,4.280)--(11.671,4.280);
\gpcolor{rgb color={0.400,0.400,0.400}}
\gpsetlinetype{gp lt border}
\gpsetdashtype{gp dt solid}
\gpsetlinewidth{2.00}
\draw[gp path] (1.839,4.280)--(1.772,4.280);
\gpcolor{rgb color={0.702,0.702,0.702}}
\gpsetlinetype{gp lt axes}
\gpsetdashtype{gp dt axes}
\gpsetlinewidth{1.00}
\draw[gp path] (1.839,4.693)--(11.671,4.693);
\gpcolor{rgb color={0.400,0.400,0.400}}
\gpsetlinetype{gp lt border}
\gpsetdashtype{gp dt solid}
\gpsetlinewidth{2.00}
\draw[gp path] (1.839,4.693)--(1.772,4.693);
\gpcolor{rgb color={0.702,0.702,0.702}}
\gpsetlinetype{gp lt axes}
\gpsetdashtype{gp dt axes}
\gpsetlinewidth{1.00}
\draw[gp path] (1.839,5.006)--(11.671,5.006);
\gpcolor{rgb color={0.400,0.400,0.400}}
\gpsetlinetype{gp lt border}
\gpsetdashtype{gp dt solid}
\gpsetlinewidth{2.00}
\draw[gp path] (1.839,5.006)--(1.772,5.006);
\gpcolor{rgb color={0.702,0.702,0.702}}
\gpsetlinetype{gp lt axes}
\gpsetdashtype{gp dt axes}
\gpsetlinewidth{1.00}
\draw[gp path] (1.839,5.259)--(11.671,5.259);
\gpcolor{rgb color={0.400,0.400,0.400}}
\gpsetlinetype{gp lt border}
\gpsetdashtype{gp dt solid}
\gpsetlinewidth{2.00}
\draw[gp path] (1.839,5.259)--(1.704,5.259);
\node[gp node right,font={\fontsize{8.0pt}{9.6pt}\selectfont}] at (1.612,5.259) {8};
\gpcolor{rgb color={0.702,0.702,0.702}}
\gpsetlinetype{gp lt axes}
\gpsetdashtype{gp dt axes}
\gpsetlinewidth{1.00}
\draw[gp path] (1.839,6.410)--(2.115,6.410);
\draw[gp path] (5.055,6.410)--(11.671,6.410);
\gpcolor{rgb color={0.400,0.400,0.400}}
\gpsetlinetype{gp lt border}
\gpsetdashtype{gp dt solid}
\gpsetlinewidth{2.00}
\draw[gp path] (1.839,6.410)--(1.772,6.410);
\gpcolor{rgb color={0.702,0.702,0.702}}
\gpsetlinetype{gp lt axes}
\gpsetdashtype{gp dt axes}
\gpsetlinewidth{1.00}
\draw[gp path] (1.839,7.015)--(2.115,7.015);
\draw[gp path] (5.055,7.015)--(11.671,7.015);
\gpcolor{rgb color={0.400,0.400,0.400}}
\gpsetlinetype{gp lt border}
\gpsetdashtype{gp dt solid}
\gpsetlinewidth{2.00}
\draw[gp path] (1.839,7.015)--(1.772,7.015);
\gpcolor{rgb color={0.702,0.702,0.702}}
\gpsetlinetype{gp lt axes}
\gpsetdashtype{gp dt axes}
\gpsetlinewidth{1.00}
\draw[gp path] (1.839,7.427)--(2.115,7.427);
\draw[gp path] (5.055,7.427)--(11.671,7.427);
\gpcolor{rgb color={0.400,0.400,0.400}}
\gpsetlinetype{gp lt border}
\gpsetdashtype{gp dt solid}
\gpsetlinewidth{2.00}
\draw[gp path] (1.839,7.427)--(1.772,7.427);
\gpcolor{rgb color={0.702,0.702,0.702}}
\gpsetlinetype{gp lt axes}
\gpsetdashtype{gp dt axes}
\gpsetlinewidth{1.00}
\draw[gp path] (1.839,7.741)--(2.115,7.741);
\draw[gp path] (5.055,7.741)--(11.671,7.741);
\gpcolor{rgb color={0.400,0.400,0.400}}
\gpsetlinetype{gp lt border}
\gpsetdashtype{gp dt solid}
\gpsetlinewidth{2.00}
\draw[gp path] (1.839,7.741)--(1.772,7.741);
\gpcolor{rgb color={0.702,0.702,0.702}}
\gpsetlinetype{gp lt axes}
\gpsetdashtype{gp dt axes}
\gpsetlinewidth{1.00}
\draw[gp path] (1.839,7.994)--(2.115,7.994);
\draw[gp path] (5.055,7.994)--(11.671,7.994);
\gpcolor{rgb color={0.400,0.400,0.400}}
\gpsetlinetype{gp lt border}
\gpsetdashtype{gp dt solid}
\gpsetlinewidth{2.00}
\draw[gp path] (1.839,7.994)--(1.704,7.994);
\node[gp node right,font={\fontsize{8.0pt}{9.6pt}\selectfont}] at (1.612,7.994) {64};
\gpcolor{rgb color={0.702,0.702,0.702}}
\gpsetlinetype{gp lt axes}
\gpsetdashtype{gp dt axes}
\gpsetlinewidth{1.00}
\draw[gp path] (2.422,1.613)--(2.422,6.028);
\draw[gp path] (2.422,8.107)--(2.422,8.287);
\gpcolor{rgb color={0.400,0.400,0.400}}
\gpsetlinetype{gp lt border}
\gpsetdashtype{gp dt solid}
\gpsetlinewidth{2.00}
\draw[gp path] (2.422,1.613)--(2.422,1.478);
\node[gp node center,font={\fontsize{8.0pt}{9.6pt}\selectfont}] at (2.422,1.016) {40};
\gpcolor{rgb color={0.702,0.702,0.702}}
\gpsetlinetype{gp lt axes}
\gpsetdashtype{gp dt axes}
\gpsetlinewidth{1.00}
\draw[gp path] (5.231,1.613)--(5.231,8.287);
\gpcolor{rgb color={0.400,0.400,0.400}}
\gpsetlinetype{gp lt border}
\gpsetdashtype{gp dt solid}
\gpsetlinewidth{2.00}
\draw[gp path] (5.231,1.613)--(5.231,1.478);
\node[gp node center,font={\fontsize{8.0pt}{9.6pt}\selectfont}] at (5.231,1.016) {160};
\gpcolor{rgb color={0.702,0.702,0.702}}
\gpsetlinetype{gp lt axes}
\gpsetdashtype{gp dt axes}
\gpsetlinewidth{1.00}
\draw[gp path] (8.040,1.613)--(8.040,8.287);
\gpcolor{rgb color={0.400,0.400,0.400}}
\gpsetlinetype{gp lt border}
\gpsetdashtype{gp dt solid}
\gpsetlinewidth{2.00}
\draw[gp path] (8.040,1.613)--(8.040,1.478);
\node[gp node center,font={\fontsize{8.0pt}{9.6pt}\selectfont}] at (8.040,1.016) {640};
\gpcolor{rgb color={0.702,0.702,0.702}}
\gpsetlinetype{gp lt axes}
\gpsetdashtype{gp dt axes}
\gpsetlinewidth{1.00}
\draw[gp path] (10.849,1.613)--(10.849,8.287);
\gpcolor{rgb color={0.400,0.400,0.400}}
\gpsetlinetype{gp lt border}
\gpsetdashtype{gp dt solid}
\gpsetlinewidth{2.00}
\draw[gp path] (10.849,1.613)--(10.849,1.478);
\node[gp node center,font={\fontsize{8.0pt}{9.6pt}\selectfont}] at (10.849,1.016) {2560};
\draw[gp path] (1.839,8.287)--(1.839,1.613)--(11.671,1.613)--(11.671,8.287)--cycle;
\gpcolor{color=gp lt color border}
\node[gp node center,rotate=-270] at (0.874,4.950) {Speedup (strong)};
\node[gp node center] at (6.755,0.400) {Number of processes};
\gpcolor{rgb color={0.400,0.400,0.400}}
\gpsetlinewidth{1.00}
\draw[gp path] (2.115,6.028)--(2.115,8.107)--(5.055,8.107)--(5.055,6.028)--cycle;
\gpcolor{color=gp lt color border}
\node[gp node left] at (3.123,7.760) {$S_{\text{ass}}$};
\gpcolor{rgb color={0.800,0.000,0.063}}
\gpsetlinewidth{3.00}
\draw[gp path] (2.391,7.760)--(2.847,7.760);
\draw[gp path] (2.422,2.525)--(5.231,4.157)--(8.040,5.820)--(10.849,7.546);
\gpsetpointsize{5.20}
\gppoint{gp mark 1}{(2.422,2.525)}
\gppoint{gp mark 1}{(5.231,4.157)}
\gppoint{gp mark 1}{(8.040,5.820)}
\gppoint{gp mark 1}{(10.849,7.546)}
\gppoint{gp mark 1}{(2.619,7.760)}
\gpcolor{color=gp lt color border}
\node[gp node left] at (3.123,7.298) {$S_{\text{init}}$};
\gpcolor{rgb color={0.800,0.000,0.663}}
\draw[gp path] (2.391,7.298)--(2.847,7.298);
\draw[gp path] (2.422,2.525)--(5.231,4.173)--(8.040,5.781)--(10.849,7.140);
\gppoint{gp mark 2}{(2.422,2.525)}
\gppoint{gp mark 2}{(5.231,4.173)}
\gppoint{gp mark 2}{(8.040,5.781)}
\gppoint{gp mark 2}{(10.849,7.140)}
\gppoint{gp mark 2}{(2.619,7.298)}
\gpcolor{color=gp lt color border}
\node[gp node left] at (3.123,6.836) {$S_{\text{solve}}$};
\gpcolor{rgb color={0.329,0.000,0.800}}
\draw[gp path] (2.391,6.836)--(2.847,6.836);
\draw[gp path] (2.422,2.525)--(5.231,4.386)--(8.040,6.194)--(10.849,7.656);
\gppoint{gp mark 3}{(2.422,2.525)}
\gppoint{gp mark 3}{(5.231,4.386)}
\gppoint{gp mark 3}{(8.040,6.194)}
\gppoint{gp mark 3}{(10.849,7.656)}
\gppoint{gp mark 3}{(2.619,6.836)}
\gpcolor{color=gp lt color border}
\node[gp node left] at (3.123,6.374) {$S_{\text{ideal}}$};
\gpcolor{rgb color={0.800,0.533,0.000}}
\draw[gp path] (2.391,6.374)--(2.847,6.374);
\draw[gp path] (2.422,2.525)--(5.231,4.348)--(8.040,6.171)--(10.849,7.994);
%% coordinates of the plot area
\gpdefrectangularnode{gp plot 1}{\pgfpoint{1.839cm}{1.613cm}}{\pgfpoint{11.671cm}{8.287cm}}
\end{tikzpicture}
%% gnuplot variables
}
    \caption{Linear elasticity problem}
  \end{subfigure}
  \begin{subfigure}[t]{\textwidth}
    \centering
    \scalebox{0.48}{\begin{tikzpicture}[gnuplot]
%% generated with GNUPLOT 5.2p6 (Lua 5.3; terminal rev. Nov 2018, script rev. 107)
%% Fr  8 Mär 17:47:29 2019
\tikzset{every node/.append style={scale=1.50}}
\path (0.000,0.000) rectangle (12.500,8.750);
\gpcolor{rgb color={0.702,0.702,0.702}}
\gpsetlinetype{gp lt axes}
\gpsetdashtype{gp dt axes}
\gpsetlinewidth{1.00}
\draw[gp path] (2.391,1.852)--(11.671,1.852);
\gpcolor{rgb color={0.400,0.400,0.400}}
\gpsetlinetype{gp lt border}
\gpsetdashtype{gp dt solid}
\gpsetlinewidth{2.00}
\draw[gp path] (2.391,1.852)--(2.324,1.852);
\gpcolor{rgb color={0.702,0.702,0.702}}
\gpsetlinetype{gp lt axes}
\gpsetdashtype{gp dt axes}
\gpsetlinewidth{1.00}
\draw[gp path] (2.391,2.091)--(11.671,2.091);
\gpcolor{rgb color={0.400,0.400,0.400}}
\gpsetlinetype{gp lt border}
\gpsetdashtype{gp dt solid}
\gpsetlinewidth{2.00}
\draw[gp path] (2.391,2.091)--(2.324,2.091);
\gpcolor{rgb color={0.702,0.702,0.702}}
\gpsetlinetype{gp lt axes}
\gpsetdashtype{gp dt axes}
\gpsetlinewidth{1.00}
\draw[gp path] (2.391,2.283)--(11.671,2.283);
\gpcolor{rgb color={0.400,0.400,0.400}}
\gpsetlinetype{gp lt border}
\gpsetdashtype{gp dt solid}
\gpsetlinewidth{2.00}
\draw[gp path] (2.391,2.283)--(2.256,2.283);
\node[gp node right,font={\fontsize{8.0pt}{9.6pt}\selectfont}] at (2.256,2.283) {10};
\gpcolor{rgb color={0.702,0.702,0.702}}
\gpsetlinetype{gp lt axes}
\gpsetdashtype{gp dt axes}
\gpsetlinewidth{1.00}
\draw[gp path] (2.391,3.277)--(11.671,3.277);
\gpcolor{rgb color={0.400,0.400,0.400}}
\gpsetlinetype{gp lt border}
\gpsetdashtype{gp dt solid}
\gpsetlinewidth{2.00}
\draw[gp path] (2.391,3.277)--(2.324,3.277);
\gpcolor{rgb color={0.702,0.702,0.702}}
\gpsetlinetype{gp lt axes}
\gpsetdashtype{gp dt axes}
\gpsetlinewidth{1.00}
\draw[gp path] (2.391,3.757)--(11.671,3.757);
\gpcolor{rgb color={0.400,0.400,0.400}}
\gpsetlinetype{gp lt border}
\gpsetdashtype{gp dt solid}
\gpsetlinewidth{2.00}
\draw[gp path] (2.391,3.757)--(2.324,3.757);
\gpcolor{rgb color={0.702,0.702,0.702}}
\gpsetlinetype{gp lt axes}
\gpsetdashtype{gp dt axes}
\gpsetlinewidth{1.00}
\draw[gp path] (2.391,4.076)--(11.671,4.076);
\gpcolor{rgb color={0.400,0.400,0.400}}
\gpsetlinetype{gp lt border}
\gpsetdashtype{gp dt solid}
\gpsetlinewidth{2.00}
\draw[gp path] (2.391,4.076)--(2.324,4.076);
\gpcolor{rgb color={0.702,0.702,0.702}}
\gpsetlinetype{gp lt axes}
\gpsetdashtype{gp dt axes}
\gpsetlinewidth{1.00}
\draw[gp path] (2.391,4.316)--(11.671,4.316);
\gpcolor{rgb color={0.400,0.400,0.400}}
\gpsetlinetype{gp lt border}
\gpsetdashtype{gp dt solid}
\gpsetlinewidth{2.00}
\draw[gp path] (2.391,4.316)--(2.324,4.316);
\gpcolor{rgb color={0.702,0.702,0.702}}
\gpsetlinetype{gp lt axes}
\gpsetdashtype{gp dt axes}
\gpsetlinewidth{1.00}
\draw[gp path] (2.391,4.507)--(11.671,4.507);
\gpcolor{rgb color={0.400,0.400,0.400}}
\gpsetlinetype{gp lt border}
\gpsetdashtype{gp dt solid}
\gpsetlinewidth{2.00}
\draw[gp path] (2.391,4.507)--(2.256,4.507);
\node[gp node right,font={\fontsize{8.0pt}{9.6pt}\selectfont}] at (2.164,4.507) {100};
\gpcolor{rgb color={0.702,0.702,0.702}}
\gpsetlinetype{gp lt axes}
\gpsetdashtype{gp dt axes}
\gpsetlinewidth{1.00}
\draw[gp path] (2.391,5.502)--(11.671,5.502);
\gpcolor{rgb color={0.400,0.400,0.400}}
\gpsetlinetype{gp lt border}
\gpsetdashtype{gp dt solid}
\gpsetlinewidth{2.00}
\draw[gp path] (2.391,5.502)--(2.324,5.502);
\gpcolor{rgb color={0.702,0.702,0.702}}
\gpsetlinetype{gp lt axes}
\gpsetdashtype{gp dt axes}
\gpsetlinewidth{1.00}
\draw[gp path] (2.391,5.982)--(11.671,5.982);
\gpcolor{rgb color={0.400,0.400,0.400}}
\gpsetlinetype{gp lt border}
\gpsetdashtype{gp dt solid}
\gpsetlinewidth{2.00}
\draw[gp path] (2.391,5.982)--(2.324,5.982);
\gpcolor{rgb color={0.702,0.702,0.702}}
\gpsetlinetype{gp lt axes}
\gpsetdashtype{gp dt axes}
\gpsetlinewidth{1.00}
\draw[gp path] (2.391,6.301)--(11.671,6.301);
\gpcolor{rgb color={0.400,0.400,0.400}}
\gpsetlinetype{gp lt border}
\gpsetdashtype{gp dt solid}
\gpsetlinewidth{2.00}
\draw[gp path] (2.391,6.301)--(2.324,6.301);
\gpcolor{rgb color={0.702,0.702,0.702}}
\gpsetlinetype{gp lt axes}
\gpsetdashtype{gp dt axes}
\gpsetlinewidth{1.00}
\draw[gp path] (2.391,6.540)--(7.075,6.540);
\draw[gp path] (11.395,6.540)--(11.671,6.540);
\gpcolor{rgb color={0.400,0.400,0.400}}
\gpsetlinetype{gp lt border}
\gpsetdashtype{gp dt solid}
\gpsetlinewidth{2.00}
\draw[gp path] (2.391,6.540)--(2.324,6.540);
\gpcolor{rgb color={0.702,0.702,0.702}}
\gpsetlinetype{gp lt axes}
\gpsetdashtype{gp dt axes}
\gpsetlinewidth{1.00}
\draw[gp path] (2.391,6.732)--(7.075,6.732);
\draw[gp path] (11.395,6.732)--(11.671,6.732);
\gpcolor{rgb color={0.400,0.400,0.400}}
\gpsetlinetype{gp lt border}
\gpsetdashtype{gp dt solid}
\gpsetlinewidth{2.00}
\draw[gp path] (2.391,6.732)--(2.256,6.732);
\node[gp node right,font={\fontsize{8.0pt}{9.6pt}\selectfont}] at (2.164,6.732) {1000};
\gpcolor{rgb color={0.702,0.702,0.702}}
\gpsetlinetype{gp lt axes}
\gpsetdashtype{gp dt axes}
\gpsetlinewidth{1.00}
\draw[gp path] (2.391,7.727)--(7.075,7.727);
\draw[gp path] (11.395,7.727)--(11.671,7.727);
\gpcolor{rgb color={0.400,0.400,0.400}}
\gpsetlinetype{gp lt border}
\gpsetdashtype{gp dt solid}
\gpsetlinewidth{2.00}
\draw[gp path] (2.391,7.727)--(2.324,7.727);
\gpcolor{rgb color={0.702,0.702,0.702}}
\gpsetlinetype{gp lt axes}
\gpsetdashtype{gp dt axes}
\gpsetlinewidth{1.00}
\draw[gp path] (2.391,8.206)--(11.671,8.206);
\gpcolor{rgb color={0.400,0.400,0.400}}
\gpsetlinetype{gp lt border}
\gpsetdashtype{gp dt solid}
\gpsetlinewidth{2.00}
\draw[gp path] (2.391,8.206)--(2.324,8.206);
\gpcolor{rgb color={0.702,0.702,0.702}}
\gpsetlinetype{gp lt axes}
\gpsetdashtype{gp dt axes}
\gpsetlinewidth{1.00}
\draw[gp path] (2.941,1.613)--(2.941,8.287);
\gpcolor{rgb color={0.400,0.400,0.400}}
\gpsetlinetype{gp lt border}
\gpsetdashtype{gp dt solid}
\gpsetlinewidth{2.00}
\draw[gp path] (2.941,1.613)--(2.941,1.478);
\node[gp node center,font={\fontsize{8.0pt}{9.6pt}\selectfont}] at (2.941,1.016) {40};
\gpcolor{rgb color={0.702,0.702,0.702}}
\gpsetlinetype{gp lt axes}
\gpsetdashtype{gp dt axes}
\gpsetlinewidth{1.00}
\draw[gp path] (5.593,1.613)--(5.593,8.287);
\gpcolor{rgb color={0.400,0.400,0.400}}
\gpsetlinetype{gp lt border}
\gpsetdashtype{gp dt solid}
\gpsetlinewidth{2.00}
\draw[gp path] (5.593,1.613)--(5.593,1.478);
\node[gp node center,font={\fontsize{8.0pt}{9.6pt}\selectfont}] at (5.593,1.016) {160};
\gpcolor{rgb color={0.702,0.702,0.702}}
\gpsetlinetype{gp lt axes}
\gpsetdashtype{gp dt axes}
\gpsetlinewidth{1.00}
\draw[gp path] (8.244,1.613)--(8.244,6.490);
\draw[gp path] (8.244,8.107)--(8.244,8.287);
\gpcolor{rgb color={0.400,0.400,0.400}}
\gpsetlinetype{gp lt border}
\gpsetdashtype{gp dt solid}
\gpsetlinewidth{2.00}
\draw[gp path] (8.244,1.613)--(8.244,1.478);
\node[gp node center,font={\fontsize{8.0pt}{9.6pt}\selectfont}] at (8.244,1.016) {640};
\gpcolor{rgb color={0.702,0.702,0.702}}
\gpsetlinetype{gp lt axes}
\gpsetdashtype{gp dt axes}
\gpsetlinewidth{1.00}
\draw[gp path] (10.896,1.613)--(10.896,6.490);
\draw[gp path] (10.896,8.107)--(10.896,8.287);
\gpcolor{rgb color={0.400,0.400,0.400}}
\gpsetlinetype{gp lt border}
\gpsetdashtype{gp dt solid}
\gpsetlinewidth{2.00}
\draw[gp path] (10.896,1.613)--(10.896,1.478);
\node[gp node center,font={\fontsize{8.0pt}{9.6pt}\selectfont}] at (10.896,1.016) {2560};
\draw[gp path] (2.391,8.287)--(2.391,1.613)--(11.671,1.613)--(11.671,8.287)--cycle;
\gpcolor{color=gp lt color border}
\node[gp node center,rotate=-270] at (0.874,4.950) {Time [s]};
\node[gp node center] at (7.031,0.400) {Number of processes};
\gpcolor{rgb color={0.400,0.400,0.400}}
\gpsetlinewidth{1.00}
\draw[gp path] (7.075,6.490)--(7.075,8.107)--(11.395,8.107)--(11.395,6.490)--cycle;
\gpcolor{color=gp lt color border}
\node[gp node left] at (8.083,7.760) {$T_{\text{lin,init}}$};
\gpcolor{rgb color={0.800,0.000,0.063}}
\gpsetlinewidth{3.00}
\draw[gp path] (7.351,7.760)--(7.807,7.760);
\draw[gp path] (2.941,5.838)--(5.593,4.590)--(8.244,3.434)--(10.896,2.546);
\gpsetpointsize{5.20}
\gppoint{gp mark 1}{(2.941,5.838)}
\gppoint{gp mark 1}{(5.593,4.590)}
\gppoint{gp mark 1}{(8.244,3.434)}
\gppoint{gp mark 1}{(10.896,2.546)}
\gppoint{gp mark 1}{(7.579,7.760)}
\gpcolor{color=gp lt color border}
\node[gp node left] at (8.083,7.298) {$T_{\text{lin,solve}}$};
\gpcolor{rgb color={0.800,0.000,0.663}}
\draw[gp path] (7.351,7.298)--(7.807,7.298);
\draw[gp path] (2.941,5.846)--(5.593,4.464)--(8.244,3.171)--(10.896,2.138);
\gppoint{gp mark 2}{(2.941,5.846)}
\gppoint{gp mark 2}{(5.593,4.464)}
\gppoint{gp mark 2}{(8.244,3.171)}
\gppoint{gp mark 2}{(10.896,2.138)}
\gppoint{gp mark 2}{(7.579,7.298)}
\gpcolor{color=gp lt color border}
\node[gp node left] at (8.083,6.836) {$T_{\text{newton}}$};
\gpcolor{rgb color={0.329,0.000,0.800}}
\draw[gp path] (7.351,6.836)--(7.807,6.836);
\draw[gp path] (2.941,7.404)--(5.593,6.072)--(8.244,4.795)--(10.896,3.660);
\gppoint{gp mark 3}{(2.941,7.404)}
\gppoint{gp mark 3}{(5.593,6.072)}
\gppoint{gp mark 3}{(8.244,4.795)}
\gppoint{gp mark 3}{(10.896,3.660)}
\gppoint{gp mark 3}{(7.579,6.836)}
%% coordinates of the plot area
\gpdefrectangularnode{gp plot 1}{\pgfpoint{2.391cm}{1.613cm}}{\pgfpoint{11.671cm}{8.287cm}}
\end{tikzpicture}
%% gnuplot variables
}
    \scalebox{0.48}{\begin{tikzpicture}[gnuplot]
%% generated with GNUPLOT 5.2p6 (Lua 5.3; terminal rev. Nov 2018, script rev. 107)
%% Fr  8 Mär 17:47:29 2019
\tikzset{every node/.append style={scale=1.50}}
\path (0.000,0.000) rectangle (12.500,8.750);
\gpcolor{rgb color={0.702,0.702,0.702}}
\gpsetlinetype{gp lt axes}
\gpsetdashtype{gp dt axes}
\gpsetlinewidth{1.00}
\draw[gp path] (1.839,2.525)--(11.671,2.525);
\gpcolor{rgb color={0.400,0.400,0.400}}
\gpsetlinetype{gp lt border}
\gpsetdashtype{gp dt solid}
\gpsetlinewidth{2.00}
\draw[gp path] (1.839,2.525)--(1.704,2.525);
\node[gp node right,font={\fontsize{8.0pt}{9.6pt}\selectfont}] at (1.704,2.525) {1};
\gpcolor{rgb color={0.702,0.702,0.702}}
\gpsetlinetype{gp lt axes}
\gpsetdashtype{gp dt axes}
\gpsetlinewidth{1.00}
\draw[gp path] (1.839,3.676)--(11.671,3.676);
\gpcolor{rgb color={0.400,0.400,0.400}}
\gpsetlinetype{gp lt border}
\gpsetdashtype{gp dt solid}
\gpsetlinewidth{2.00}
\draw[gp path] (1.839,3.676)--(1.772,3.676);
\gpcolor{rgb color={0.702,0.702,0.702}}
\gpsetlinetype{gp lt axes}
\gpsetdashtype{gp dt axes}
\gpsetlinewidth{1.00}
\draw[gp path] (1.839,4.280)--(11.671,4.280);
\gpcolor{rgb color={0.400,0.400,0.400}}
\gpsetlinetype{gp lt border}
\gpsetdashtype{gp dt solid}
\gpsetlinewidth{2.00}
\draw[gp path] (1.839,4.280)--(1.772,4.280);
\gpcolor{rgb color={0.702,0.702,0.702}}
\gpsetlinetype{gp lt axes}
\gpsetdashtype{gp dt axes}
\gpsetlinewidth{1.00}
\draw[gp path] (1.839,4.693)--(11.671,4.693);
\gpcolor{rgb color={0.400,0.400,0.400}}
\gpsetlinetype{gp lt border}
\gpsetdashtype{gp dt solid}
\gpsetlinewidth{2.00}
\draw[gp path] (1.839,4.693)--(1.772,4.693);
\gpcolor{rgb color={0.702,0.702,0.702}}
\gpsetlinetype{gp lt axes}
\gpsetdashtype{gp dt axes}
\gpsetlinewidth{1.00}
\draw[gp path] (1.839,5.006)--(11.671,5.006);
\gpcolor{rgb color={0.400,0.400,0.400}}
\gpsetlinetype{gp lt border}
\gpsetdashtype{gp dt solid}
\gpsetlinewidth{2.00}
\draw[gp path] (1.839,5.006)--(1.772,5.006);
\gpcolor{rgb color={0.702,0.702,0.702}}
\gpsetlinetype{gp lt axes}
\gpsetdashtype{gp dt axes}
\gpsetlinewidth{1.00}
\draw[gp path] (1.839,5.259)--(11.671,5.259);
\gpcolor{rgb color={0.400,0.400,0.400}}
\gpsetlinetype{gp lt border}
\gpsetdashtype{gp dt solid}
\gpsetlinewidth{2.00}
\draw[gp path] (1.839,5.259)--(1.704,5.259);
\node[gp node right,font={\fontsize{8.0pt}{9.6pt}\selectfont}] at (1.612,5.259) {8};
\gpcolor{rgb color={0.702,0.702,0.702}}
\gpsetlinetype{gp lt axes}
\gpsetdashtype{gp dt axes}
\gpsetlinewidth{1.00}
\draw[gp path] (1.839,6.410)--(2.115,6.410);
\draw[gp path] (6.159,6.410)--(11.671,6.410);
\gpcolor{rgb color={0.400,0.400,0.400}}
\gpsetlinetype{gp lt border}
\gpsetdashtype{gp dt solid}
\gpsetlinewidth{2.00}
\draw[gp path] (1.839,6.410)--(1.772,6.410);
\gpcolor{rgb color={0.702,0.702,0.702}}
\gpsetlinetype{gp lt axes}
\gpsetdashtype{gp dt axes}
\gpsetlinewidth{1.00}
\draw[gp path] (1.839,7.015)--(2.115,7.015);
\draw[gp path] (6.159,7.015)--(11.671,7.015);
\gpcolor{rgb color={0.400,0.400,0.400}}
\gpsetlinetype{gp lt border}
\gpsetdashtype{gp dt solid}
\gpsetlinewidth{2.00}
\draw[gp path] (1.839,7.015)--(1.772,7.015);
\gpcolor{rgb color={0.702,0.702,0.702}}
\gpsetlinetype{gp lt axes}
\gpsetdashtype{gp dt axes}
\gpsetlinewidth{1.00}
\draw[gp path] (1.839,7.427)--(2.115,7.427);
\draw[gp path] (6.159,7.427)--(11.671,7.427);
\gpcolor{rgb color={0.400,0.400,0.400}}
\gpsetlinetype{gp lt border}
\gpsetdashtype{gp dt solid}
\gpsetlinewidth{2.00}
\draw[gp path] (1.839,7.427)--(1.772,7.427);
\gpcolor{rgb color={0.702,0.702,0.702}}
\gpsetlinetype{gp lt axes}
\gpsetdashtype{gp dt axes}
\gpsetlinewidth{1.00}
\draw[gp path] (1.839,7.741)--(2.115,7.741);
\draw[gp path] (6.159,7.741)--(11.671,7.741);
\gpcolor{rgb color={0.400,0.400,0.400}}
\gpsetlinetype{gp lt border}
\gpsetdashtype{gp dt solid}
\gpsetlinewidth{2.00}
\draw[gp path] (1.839,7.741)--(1.772,7.741);
\gpcolor{rgb color={0.702,0.702,0.702}}
\gpsetlinetype{gp lt axes}
\gpsetdashtype{gp dt axes}
\gpsetlinewidth{1.00}
\draw[gp path] (1.839,7.994)--(2.115,7.994);
\draw[gp path] (6.159,7.994)--(11.671,7.994);
\gpcolor{rgb color={0.400,0.400,0.400}}
\gpsetlinetype{gp lt border}
\gpsetdashtype{gp dt solid}
\gpsetlinewidth{2.00}
\draw[gp path] (1.839,7.994)--(1.704,7.994);
\node[gp node right,font={\fontsize{8.0pt}{9.6pt}\selectfont}] at (1.612,7.994) {64};
\gpcolor{rgb color={0.702,0.702,0.702}}
\gpsetlinetype{gp lt axes}
\gpsetdashtype{gp dt axes}
\gpsetlinewidth{1.00}
\draw[gp path] (2.422,1.613)--(2.422,6.028);
\draw[gp path] (2.422,8.107)--(2.422,8.287);
\gpcolor{rgb color={0.400,0.400,0.400}}
\gpsetlinetype{gp lt border}
\gpsetdashtype{gp dt solid}
\gpsetlinewidth{2.00}
\draw[gp path] (2.422,1.613)--(2.422,1.478);
\node[gp node center,font={\fontsize{8.0pt}{9.6pt}\selectfont}] at (2.422,1.016) {40};
\gpcolor{rgb color={0.702,0.702,0.702}}
\gpsetlinetype{gp lt axes}
\gpsetdashtype{gp dt axes}
\gpsetlinewidth{1.00}
\draw[gp path] (5.231,1.613)--(5.231,6.028);
\draw[gp path] (5.231,8.107)--(5.231,8.287);
\gpcolor{rgb color={0.400,0.400,0.400}}
\gpsetlinetype{gp lt border}
\gpsetdashtype{gp dt solid}
\gpsetlinewidth{2.00}
\draw[gp path] (5.231,1.613)--(5.231,1.478);
\node[gp node center,font={\fontsize{8.0pt}{9.6pt}\selectfont}] at (5.231,1.016) {160};
\gpcolor{rgb color={0.702,0.702,0.702}}
\gpsetlinetype{gp lt axes}
\gpsetdashtype{gp dt axes}
\gpsetlinewidth{1.00}
\draw[gp path] (8.040,1.613)--(8.040,8.287);
\gpcolor{rgb color={0.400,0.400,0.400}}
\gpsetlinetype{gp lt border}
\gpsetdashtype{gp dt solid}
\gpsetlinewidth{2.00}
\draw[gp path] (8.040,1.613)--(8.040,1.478);
\node[gp node center,font={\fontsize{8.0pt}{9.6pt}\selectfont}] at (8.040,1.016) {640};
\gpcolor{rgb color={0.702,0.702,0.702}}
\gpsetlinetype{gp lt axes}
\gpsetdashtype{gp dt axes}
\gpsetlinewidth{1.00}
\draw[gp path] (10.849,1.613)--(10.849,8.287);
\gpcolor{rgb color={0.400,0.400,0.400}}
\gpsetlinetype{gp lt border}
\gpsetdashtype{gp dt solid}
\gpsetlinewidth{2.00}
\draw[gp path] (10.849,1.613)--(10.849,1.478);
\node[gp node center,font={\fontsize{8.0pt}{9.6pt}\selectfont}] at (10.849,1.016) {2560};
\draw[gp path] (1.839,8.287)--(1.839,1.613)--(11.671,1.613)--(11.671,8.287)--cycle;
\gpcolor{color=gp lt color border}
\node[gp node center,rotate=-270] at (0.874,4.950) {Speedup (strong)};
\node[gp node center] at (6.755,0.400) {Number of processes};
\gpcolor{rgb color={0.400,0.400,0.400}}
\gpsetlinewidth{1.00}
\draw[gp path] (2.115,6.028)--(2.115,8.107)--(6.159,8.107)--(6.159,6.028)--cycle;
\gpcolor{color=gp lt color border}
\node[gp node left] at (3.123,7.760) {$S_{\text{lin,init}}$};
\gpcolor{rgb color={0.800,0.000,0.063}}
\gpsetlinewidth{3.00}
\draw[gp path] (2.391,7.760)--(2.847,7.760);
\draw[gp path] (2.422,2.525)--(5.231,4.222)--(8.040,5.797)--(10.849,7.005);
\gpsetpointsize{5.20}
\gppoint{gp mark 1}{(2.422,2.525)}
\gppoint{gp mark 1}{(5.231,4.222)}
\gppoint{gp mark 1}{(8.040,5.797)}
\gppoint{gp mark 1}{(10.849,7.005)}
\gppoint{gp mark 1}{(2.619,7.760)}
\gpcolor{color=gp lt color border}
\node[gp node left] at (3.123,7.298) {$S_{\text{lin,solve}}$};
\gpcolor{rgb color={0.800,0.000,0.663}}
\draw[gp path] (2.391,7.298)--(2.847,7.298);
\draw[gp path] (2.422,2.525)--(5.231,4.406)--(8.040,6.165)--(10.849,7.571);
\gppoint{gp mark 2}{(2.422,2.525)}
\gppoint{gp mark 2}{(5.231,4.406)}
\gppoint{gp mark 2}{(8.040,6.165)}
\gppoint{gp mark 2}{(10.849,7.571)}
\gppoint{gp mark 2}{(2.619,7.298)}
\gpcolor{color=gp lt color border}
\node[gp node left] at (3.123,6.836) {$S_{\text{newton}}$};
\gpcolor{rgb color={0.329,0.000,0.800}}
\draw[gp path] (2.391,6.836)--(2.847,6.836);
\draw[gp path] (2.422,2.525)--(5.231,4.337)--(8.040,6.075)--(10.849,7.619);
\gppoint{gp mark 3}{(2.422,2.525)}
\gppoint{gp mark 3}{(5.231,4.337)}
\gppoint{gp mark 3}{(8.040,6.075)}
\gppoint{gp mark 3}{(10.849,7.619)}
\gppoint{gp mark 3}{(2.619,6.836)}
\gpcolor{color=gp lt color border}
\node[gp node left] at (3.123,6.374) {$S_{\text{ideal}}$};
\gpcolor{rgb color={0.800,0.533,0.000}}
\draw[gp path] (2.391,6.374)--(2.847,6.374);
\draw[gp path] (2.422,2.525)--(5.231,4.348)--(8.040,6.171)--(10.849,7.994);
%% coordinates of the plot area
\gpdefrectangularnode{gp plot 1}{\pgfpoint{1.839cm}{1.613cm}}{\pgfpoint{11.671cm}{8.287cm}}
\end{tikzpicture}
%% gnuplot variables
}
    \caption{Computation of shape derivative}
  \end{subfigure}
  \begin{subfigure}[t]{\textwidth}
    \centering
    \begin{tabular}{rrrccc}
    \\
    \toprule
    Procs  & Refs  & DoFs  & Linear solver  & Newton solver   & Linear solver     \\
          &       &        &  (elasticity)      & (shape derivative)    & (shape derivative)\\ 
    \midrule
    40 & 7 & 37,153,922 & 305 & 140 & 567  \\
    160 & 7 & 37,153,922 & 297 & 140 & 567  \\
    640 & 7 & 37,153,922 & 308 & 140 & 627  \\
    2560 & 7 & 37,153,922 & 306 & 140 & 630  \\
    \bottomrule
     \end{tabular}
    \caption{Accumulated iteration counts for the strong scaling study}
    \label{tab:StrongScaling}
  \end{subfigure}

    \caption{Strong scaling: accumulated wallclock time for the first 10 optimization steps (left), speedup relative to 40 cores (right), and accumulated iteration counts (table). Shown are the fine level matrix assembling (\texttt{ass}), multigrid initialization (\texttt{init}) and solution (\texttt{solve}) phase, and the overall newton algorithm (\texttt{newton})}
    \label{fig:StrongScaling}
 \end{figure}

\begin{figure}[h!tbp]

 \centering
  \begin{subfigure}[t]{\textwidth}
    \scalebox{0.49}{\input{Parmetis-NumSteps10-weakAvg-coarse-timings}} 
    \scalebox{0.49}{\begin{tikzpicture}[gnuplot]
%% generated with GNUPLOT 5.2p6 (Lua 5.3; terminal rev. Nov 2018, script rev. 107)
%% Fr  8 Mär 17:47:29 2019
\tikzset{every node/.append style={scale=1.50}}
\path (0.000,0.000) rectangle (12.500,8.750);
\gpcolor{rgb color={0.702,0.702,0.702}}
\gpsetlinetype{gp lt axes}
\gpsetdashtype{gp dt axes}
\gpsetlinewidth{1.00}
\draw[gp path] (1.839,2.525)--(11.671,2.525);
\gpcolor{rgb color={0.400,0.400,0.400}}
\gpsetlinetype{gp lt border}
\gpsetdashtype{gp dt solid}
\gpsetlinewidth{2.00}
\draw[gp path] (1.839,2.525)--(1.704,2.525);
\node[gp node right,font={\fontsize{8.0pt}{9.6pt}\selectfont}] at (1.704,2.525) {1};
\gpcolor{rgb color={0.702,0.702,0.702}}
\gpsetlinetype{gp lt axes}
\gpsetdashtype{gp dt axes}
\gpsetlinewidth{1.00}
\draw[gp path] (1.839,3.676)--(11.671,3.676);
\gpcolor{rgb color={0.400,0.400,0.400}}
\gpsetlinetype{gp lt border}
\gpsetdashtype{gp dt solid}
\gpsetlinewidth{2.00}
\draw[gp path] (1.839,3.676)--(1.772,3.676);
\gpcolor{rgb color={0.702,0.702,0.702}}
\gpsetlinetype{gp lt axes}
\gpsetdashtype{gp dt axes}
\gpsetlinewidth{1.00}
\draw[gp path] (1.839,4.280)--(11.671,4.280);
\gpcolor{rgb color={0.400,0.400,0.400}}
\gpsetlinetype{gp lt border}
\gpsetdashtype{gp dt solid}
\gpsetlinewidth{2.00}
\draw[gp path] (1.839,4.280)--(1.772,4.280);
\gpcolor{rgb color={0.702,0.702,0.702}}
\gpsetlinetype{gp lt axes}
\gpsetdashtype{gp dt axes}
\gpsetlinewidth{1.00}
\draw[gp path] (1.839,4.693)--(11.671,4.693);
\gpcolor{rgb color={0.400,0.400,0.400}}
\gpsetlinetype{gp lt border}
\gpsetdashtype{gp dt solid}
\gpsetlinewidth{2.00}
\draw[gp path] (1.839,4.693)--(1.772,4.693);
\gpcolor{rgb color={0.702,0.702,0.702}}
\gpsetlinetype{gp lt axes}
\gpsetdashtype{gp dt axes}
\gpsetlinewidth{1.00}
\draw[gp path] (1.839,5.006)--(11.671,5.006);
\gpcolor{rgb color={0.400,0.400,0.400}}
\gpsetlinetype{gp lt border}
\gpsetdashtype{gp dt solid}
\gpsetlinewidth{2.00}
\draw[gp path] (1.839,5.006)--(1.772,5.006);
\gpcolor{rgb color={0.702,0.702,0.702}}
\gpsetlinetype{gp lt axes}
\gpsetdashtype{gp dt axes}
\gpsetlinewidth{1.00}
\draw[gp path] (1.839,5.259)--(11.671,5.259);
\gpcolor{rgb color={0.400,0.400,0.400}}
\gpsetlinetype{gp lt border}
\gpsetdashtype{gp dt solid}
\gpsetlinewidth{2.00}
\draw[gp path] (1.839,5.259)--(1.704,5.259);
\node[gp node right,font={\fontsize{8.0pt}{9.6pt}\selectfont}] at (1.612,5.259) {8};
\gpcolor{rgb color={0.702,0.702,0.702}}
\gpsetlinetype{gp lt axes}
\gpsetdashtype{gp dt axes}
\gpsetlinewidth{1.00}
\draw[gp path] (1.839,6.410)--(2.115,6.410);
\draw[gp path] (5.055,6.410)--(11.671,6.410);
\gpcolor{rgb color={0.400,0.400,0.400}}
\gpsetlinetype{gp lt border}
\gpsetdashtype{gp dt solid}
\gpsetlinewidth{2.00}
\draw[gp path] (1.839,6.410)--(1.772,6.410);
\gpcolor{rgb color={0.702,0.702,0.702}}
\gpsetlinetype{gp lt axes}
\gpsetdashtype{gp dt axes}
\gpsetlinewidth{1.00}
\draw[gp path] (1.839,7.015)--(2.115,7.015);
\draw[gp path] (5.055,7.015)--(11.671,7.015);
\gpcolor{rgb color={0.400,0.400,0.400}}
\gpsetlinetype{gp lt border}
\gpsetdashtype{gp dt solid}
\gpsetlinewidth{2.00}
\draw[gp path] (1.839,7.015)--(1.772,7.015);
\gpcolor{rgb color={0.702,0.702,0.702}}
\gpsetlinetype{gp lt axes}
\gpsetdashtype{gp dt axes}
\gpsetlinewidth{1.00}
\draw[gp path] (1.839,7.427)--(2.115,7.427);
\draw[gp path] (5.055,7.427)--(11.671,7.427);
\gpcolor{rgb color={0.400,0.400,0.400}}
\gpsetlinetype{gp lt border}
\gpsetdashtype{gp dt solid}
\gpsetlinewidth{2.00}
\draw[gp path] (1.839,7.427)--(1.772,7.427);
\gpcolor{rgb color={0.702,0.702,0.702}}
\gpsetlinetype{gp lt axes}
\gpsetdashtype{gp dt axes}
\gpsetlinewidth{1.00}
\draw[gp path] (1.839,7.741)--(2.115,7.741);
\draw[gp path] (5.055,7.741)--(11.671,7.741);
\gpcolor{rgb color={0.400,0.400,0.400}}
\gpsetlinetype{gp lt border}
\gpsetdashtype{gp dt solid}
\gpsetlinewidth{2.00}
\draw[gp path] (1.839,7.741)--(1.772,7.741);
\gpcolor{rgb color={0.702,0.702,0.702}}
\gpsetlinetype{gp lt axes}
\gpsetdashtype{gp dt axes}
\gpsetlinewidth{1.00}
\draw[gp path] (1.839,7.994)--(2.115,7.994);
\draw[gp path] (5.055,7.994)--(11.671,7.994);
\gpcolor{rgb color={0.400,0.400,0.400}}
\gpsetlinetype{gp lt border}
\gpsetdashtype{gp dt solid}
\gpsetlinewidth{2.00}
\draw[gp path] (1.839,7.994)--(1.704,7.994);
\node[gp node right,font={\fontsize{8.0pt}{9.6pt}\selectfont}] at (1.612,7.994) {64};
\gpcolor{rgb color={0.702,0.702,0.702}}
\gpsetlinetype{gp lt axes}
\gpsetdashtype{gp dt axes}
\gpsetlinewidth{1.00}
\draw[gp path] (2.422,1.613)--(2.422,6.028);
\draw[gp path] (2.422,8.107)--(2.422,8.287);
\gpcolor{rgb color={0.400,0.400,0.400}}
\gpsetlinetype{gp lt border}
\gpsetdashtype{gp dt solid}
\gpsetlinewidth{2.00}
\draw[gp path] (2.422,1.613)--(2.422,1.478);
\node[gp node center,font={\fontsize{8.0pt}{9.6pt}\selectfont}] at (2.422,1.016) {40};
\gpcolor{rgb color={0.702,0.702,0.702}}
\gpsetlinetype{gp lt axes}
\gpsetdashtype{gp dt axes}
\gpsetlinewidth{1.00}
\draw[gp path] (5.231,1.613)--(5.231,8.287);
\gpcolor{rgb color={0.400,0.400,0.400}}
\gpsetlinetype{gp lt border}
\gpsetdashtype{gp dt solid}
\gpsetlinewidth{2.00}
\draw[gp path] (5.231,1.613)--(5.231,1.478);
\node[gp node center,font={\fontsize{8.0pt}{9.6pt}\selectfont}] at (5.231,1.016) {160};
\gpcolor{rgb color={0.702,0.702,0.702}}
\gpsetlinetype{gp lt axes}
\gpsetdashtype{gp dt axes}
\gpsetlinewidth{1.00}
\draw[gp path] (8.040,1.613)--(8.040,8.287);
\gpcolor{rgb color={0.400,0.400,0.400}}
\gpsetlinetype{gp lt border}
\gpsetdashtype{gp dt solid}
\gpsetlinewidth{2.00}
\draw[gp path] (8.040,1.613)--(8.040,1.478);
\node[gp node center,font={\fontsize{8.0pt}{9.6pt}\selectfont}] at (8.040,1.016) {640};
\gpcolor{rgb color={0.702,0.702,0.702}}
\gpsetlinetype{gp lt axes}
\gpsetdashtype{gp dt axes}
\gpsetlinewidth{1.00}
\draw[gp path] (10.849,1.613)--(10.849,8.287);
\gpcolor{rgb color={0.400,0.400,0.400}}
\gpsetlinetype{gp lt border}
\gpsetdashtype{gp dt solid}
\gpsetlinewidth{2.00}
\draw[gp path] (10.849,1.613)--(10.849,1.478);
\node[gp node center,font={\fontsize{8.0pt}{9.6pt}\selectfont}] at (10.849,1.016) {2560};
\draw[gp path] (1.839,8.287)--(1.839,1.613)--(11.671,1.613)--(11.671,8.287)--cycle;
\gpcolor{color=gp lt color border}
\node[gp node center,rotate=-270] at (0.874,4.950) {Speedup (weak)};
\node[gp node center] at (6.755,0.400) {Number of processes};
\gpcolor{rgb color={0.400,0.400,0.400}}
\gpsetlinewidth{1.00}
\draw[gp path] (2.115,6.028)--(2.115,8.107)--(5.055,8.107)--(5.055,6.028)--cycle;
\gpcolor{color=gp lt color border}
\node[gp node left] at (3.123,7.760) {$S_{\text{ass}}$};
\gpcolor{rgb color={0.800,0.000,0.063}}
\gpsetlinewidth{3.00}
\draw[gp path] (2.391,7.760)--(2.847,7.760);
\draw[gp path] (2.422,2.525)--(5.231,4.173)--(8.040,5.823)--(10.849,7.567);
\gpsetpointsize{5.20}
\gppoint{gp mark 1}{(2.422,2.525)}
\gppoint{gp mark 1}{(5.231,4.173)}
\gppoint{gp mark 1}{(8.040,5.823)}
\gppoint{gp mark 1}{(10.849,7.567)}
\gppoint{gp mark 1}{(2.619,7.760)}
\gpcolor{color=gp lt color border}
\node[gp node left] at (3.123,7.298) {$S_{\text{init}}$};
\gpcolor{rgb color={0.800,0.000,0.663}}
\draw[gp path] (2.391,7.298)--(2.847,7.298);
\draw[gp path] (2.422,2.525)--(5.231,4.202)--(8.040,5.880)--(10.849,7.620);
\gppoint{gp mark 2}{(2.422,2.525)}
\gppoint{gp mark 2}{(5.231,4.202)}
\gppoint{gp mark 2}{(8.040,5.880)}
\gppoint{gp mark 2}{(10.849,7.620)}
\gppoint{gp mark 2}{(2.619,7.298)}
\gpcolor{color=gp lt color border}
\node[gp node left] at (3.123,6.836) {$S_{\text{solve}}$};
\gpcolor{rgb color={0.329,0.000,0.800}}
\draw[gp path] (2.391,6.836)--(2.847,6.836);
\draw[gp path] (2.422,2.525)--(5.231,4.299)--(8.040,6.061)--(10.849,7.796);
\gppoint{gp mark 3}{(2.422,2.525)}
\gppoint{gp mark 3}{(5.231,4.299)}
\gppoint{gp mark 3}{(8.040,6.061)}
\gppoint{gp mark 3}{(10.849,7.796)}
\gppoint{gp mark 3}{(2.619,6.836)}
\gpcolor{color=gp lt color border}
\node[gp node left] at (3.123,6.374) {$S_{\text{ideal}}$};
\gpcolor{rgb color={0.800,0.533,0.000}}
\draw[gp path] (2.391,6.374)--(2.847,6.374);
\draw[gp path] (2.422,2.525)--(5.231,4.348)--(8.040,6.171)--(10.849,7.994);
%% coordinates of the plot area
\gpdefrectangularnode{gp plot 1}{\pgfpoint{1.839cm}{1.613cm}}{\pgfpoint{11.671cm}{8.287cm}}
\end{tikzpicture}
%% gnuplot variables
}
    \caption{Linear elasticity problem}
  \end{subfigure}
 \begin{subfigure}[t]{\textwidth}
    \scalebox{0.49}{\input{Parmetis-NumSteps10-shape-weakAvg-coarse-timings}} 
    \scalebox{0.49}{\begin{tikzpicture}[gnuplot]
%% generated with GNUPLOT 5.2p6 (Lua 5.3; terminal rev. Nov 2018, script rev. 107)
%% Fr  8 Mär 17:47:29 2019
\tikzset{every node/.append style={scale=1.50}}
\path (0.000,0.000) rectangle (12.500,8.750);
\gpcolor{rgb color={0.702,0.702,0.702}}
\gpsetlinetype{gp lt axes}
\gpsetdashtype{gp dt axes}
\gpsetlinewidth{1.00}
\draw[gp path] (1.839,2.525)--(11.671,2.525);
\gpcolor{rgb color={0.400,0.400,0.400}}
\gpsetlinetype{gp lt border}
\gpsetdashtype{gp dt solid}
\gpsetlinewidth{2.00}
\draw[gp path] (1.839,2.525)--(1.704,2.525);
\node[gp node right,font={\fontsize{8.0pt}{9.6pt}\selectfont}] at (1.704,2.525) {1};
\gpcolor{rgb color={0.702,0.702,0.702}}
\gpsetlinetype{gp lt axes}
\gpsetdashtype{gp dt axes}
\gpsetlinewidth{1.00}
\draw[gp path] (1.839,3.676)--(11.671,3.676);
\gpcolor{rgb color={0.400,0.400,0.400}}
\gpsetlinetype{gp lt border}
\gpsetdashtype{gp dt solid}
\gpsetlinewidth{2.00}
\draw[gp path] (1.839,3.676)--(1.772,3.676);
\gpcolor{rgb color={0.702,0.702,0.702}}
\gpsetlinetype{gp lt axes}
\gpsetdashtype{gp dt axes}
\gpsetlinewidth{1.00}
\draw[gp path] (1.839,4.280)--(11.671,4.280);
\gpcolor{rgb color={0.400,0.400,0.400}}
\gpsetlinetype{gp lt border}
\gpsetdashtype{gp dt solid}
\gpsetlinewidth{2.00}
\draw[gp path] (1.839,4.280)--(1.772,4.280);
\gpcolor{rgb color={0.702,0.702,0.702}}
\gpsetlinetype{gp lt axes}
\gpsetdashtype{gp dt axes}
\gpsetlinewidth{1.00}
\draw[gp path] (1.839,4.693)--(11.671,4.693);
\gpcolor{rgb color={0.400,0.400,0.400}}
\gpsetlinetype{gp lt border}
\gpsetdashtype{gp dt solid}
\gpsetlinewidth{2.00}
\draw[gp path] (1.839,4.693)--(1.772,4.693);
\gpcolor{rgb color={0.702,0.702,0.702}}
\gpsetlinetype{gp lt axes}
\gpsetdashtype{gp dt axes}
\gpsetlinewidth{1.00}
\draw[gp path] (1.839,5.006)--(11.671,5.006);
\gpcolor{rgb color={0.400,0.400,0.400}}
\gpsetlinetype{gp lt border}
\gpsetdashtype{gp dt solid}
\gpsetlinewidth{2.00}
\draw[gp path] (1.839,5.006)--(1.772,5.006);
\gpcolor{rgb color={0.702,0.702,0.702}}
\gpsetlinetype{gp lt axes}
\gpsetdashtype{gp dt axes}
\gpsetlinewidth{1.00}
\draw[gp path] (1.839,5.259)--(11.671,5.259);
\gpcolor{rgb color={0.400,0.400,0.400}}
\gpsetlinetype{gp lt border}
\gpsetdashtype{gp dt solid}
\gpsetlinewidth{2.00}
\draw[gp path] (1.839,5.259)--(1.704,5.259);
\node[gp node right,font={\fontsize{8.0pt}{9.6pt}\selectfont}] at (1.612,5.259) {8};
\gpcolor{rgb color={0.702,0.702,0.702}}
\gpsetlinetype{gp lt axes}
\gpsetdashtype{gp dt axes}
\gpsetlinewidth{1.00}
\draw[gp path] (1.839,6.410)--(2.115,6.410);
\draw[gp path] (6.159,6.410)--(11.671,6.410);
\gpcolor{rgb color={0.400,0.400,0.400}}
\gpsetlinetype{gp lt border}
\gpsetdashtype{gp dt solid}
\gpsetlinewidth{2.00}
\draw[gp path] (1.839,6.410)--(1.772,6.410);
\gpcolor{rgb color={0.702,0.702,0.702}}
\gpsetlinetype{gp lt axes}
\gpsetdashtype{gp dt axes}
\gpsetlinewidth{1.00}
\draw[gp path] (1.839,7.015)--(2.115,7.015);
\draw[gp path] (6.159,7.015)--(11.671,7.015);
\gpcolor{rgb color={0.400,0.400,0.400}}
\gpsetlinetype{gp lt border}
\gpsetdashtype{gp dt solid}
\gpsetlinewidth{2.00}
\draw[gp path] (1.839,7.015)--(1.772,7.015);
\gpcolor{rgb color={0.702,0.702,0.702}}
\gpsetlinetype{gp lt axes}
\gpsetdashtype{gp dt axes}
\gpsetlinewidth{1.00}
\draw[gp path] (1.839,7.427)--(2.115,7.427);
\draw[gp path] (6.159,7.427)--(11.671,7.427);
\gpcolor{rgb color={0.400,0.400,0.400}}
\gpsetlinetype{gp lt border}
\gpsetdashtype{gp dt solid}
\gpsetlinewidth{2.00}
\draw[gp path] (1.839,7.427)--(1.772,7.427);
\gpcolor{rgb color={0.702,0.702,0.702}}
\gpsetlinetype{gp lt axes}
\gpsetdashtype{gp dt axes}
\gpsetlinewidth{1.00}
\draw[gp path] (1.839,7.741)--(2.115,7.741);
\draw[gp path] (6.159,7.741)--(11.671,7.741);
\gpcolor{rgb color={0.400,0.400,0.400}}
\gpsetlinetype{gp lt border}
\gpsetdashtype{gp dt solid}
\gpsetlinewidth{2.00}
\draw[gp path] (1.839,7.741)--(1.772,7.741);
\gpcolor{rgb color={0.702,0.702,0.702}}
\gpsetlinetype{gp lt axes}
\gpsetdashtype{gp dt axes}
\gpsetlinewidth{1.00}
\draw[gp path] (1.839,7.994)--(2.115,7.994);
\draw[gp path] (6.159,7.994)--(11.671,7.994);
\gpcolor{rgb color={0.400,0.400,0.400}}
\gpsetlinetype{gp lt border}
\gpsetdashtype{gp dt solid}
\gpsetlinewidth{2.00}
\draw[gp path] (1.839,7.994)--(1.704,7.994);
\node[gp node right,font={\fontsize{8.0pt}{9.6pt}\selectfont}] at (1.612,7.994) {64};
\gpcolor{rgb color={0.702,0.702,0.702}}
\gpsetlinetype{gp lt axes}
\gpsetdashtype{gp dt axes}
\gpsetlinewidth{1.00}
\draw[gp path] (2.422,1.613)--(2.422,6.028);
\draw[gp path] (2.422,8.107)--(2.422,8.287);
\gpcolor{rgb color={0.400,0.400,0.400}}
\gpsetlinetype{gp lt border}
\gpsetdashtype{gp dt solid}
\gpsetlinewidth{2.00}
\draw[gp path] (2.422,1.613)--(2.422,1.478);
\node[gp node center,font={\fontsize{8.0pt}{9.6pt}\selectfont}] at (2.422,1.016) {40};
\gpcolor{rgb color={0.702,0.702,0.702}}
\gpsetlinetype{gp lt axes}
\gpsetdashtype{gp dt axes}
\gpsetlinewidth{1.00}
\draw[gp path] (5.231,1.613)--(5.231,6.028);
\draw[gp path] (5.231,8.107)--(5.231,8.287);
\gpcolor{rgb color={0.400,0.400,0.400}}
\gpsetlinetype{gp lt border}
\gpsetdashtype{gp dt solid}
\gpsetlinewidth{2.00}
\draw[gp path] (5.231,1.613)--(5.231,1.478);
\node[gp node center,font={\fontsize{8.0pt}{9.6pt}\selectfont}] at (5.231,1.016) {160};
\gpcolor{rgb color={0.702,0.702,0.702}}
\gpsetlinetype{gp lt axes}
\gpsetdashtype{gp dt axes}
\gpsetlinewidth{1.00}
\draw[gp path] (8.040,1.613)--(8.040,8.287);
\gpcolor{rgb color={0.400,0.400,0.400}}
\gpsetlinetype{gp lt border}
\gpsetdashtype{gp dt solid}
\gpsetlinewidth{2.00}
\draw[gp path] (8.040,1.613)--(8.040,1.478);
\node[gp node center,font={\fontsize{8.0pt}{9.6pt}\selectfont}] at (8.040,1.016) {640};
\gpcolor{rgb color={0.702,0.702,0.702}}
\gpsetlinetype{gp lt axes}
\gpsetdashtype{gp dt axes}
\gpsetlinewidth{1.00}
\draw[gp path] (10.849,1.613)--(10.849,8.287);
\gpcolor{rgb color={0.400,0.400,0.400}}
\gpsetlinetype{gp lt border}
\gpsetdashtype{gp dt solid}
\gpsetlinewidth{2.00}
\draw[gp path] (10.849,1.613)--(10.849,1.478);
\node[gp node center,font={\fontsize{8.0pt}{9.6pt}\selectfont}] at (10.849,1.016) {2560};
\draw[gp path] (1.839,8.287)--(1.839,1.613)--(11.671,1.613)--(11.671,8.287)--cycle;
\gpcolor{color=gp lt color border}
\node[gp node center,rotate=-270] at (0.874,4.950) {Speedup (weak)};
\node[gp node center] at (6.755,0.400) {Number of processes};
\gpcolor{rgb color={0.400,0.400,0.400}}
\gpsetlinewidth{1.00}
\draw[gp path] (2.115,6.028)--(2.115,8.107)--(6.159,8.107)--(6.159,6.028)--cycle;
\gpcolor{color=gp lt color border}
\node[gp node left] at (3.123,7.760) {$S_{\text{lin,init}}$};
\gpcolor{rgb color={0.800,0.000,0.063}}
\gpsetlinewidth{3.00}
\draw[gp path] (2.391,7.760)--(2.847,7.760);
\draw[gp path] (2.422,2.525)--(5.231,4.193)--(8.040,5.923)--(10.849,7.618);
\gpsetpointsize{5.20}
\gppoint{gp mark 1}{(2.422,2.525)}
\gppoint{gp mark 1}{(5.231,4.193)}
\gppoint{gp mark 1}{(8.040,5.923)}
\gppoint{gp mark 1}{(10.849,7.618)}
\gppoint{gp mark 1}{(2.619,7.760)}
\gpcolor{color=gp lt color border}
\node[gp node left] at (3.123,7.298) {$S_{\text{lin,solve}}$};
\gpcolor{rgb color={0.800,0.000,0.663}}
\draw[gp path] (2.391,7.298)--(2.847,7.298);
\draw[gp path] (2.422,2.525)--(5.231,4.218)--(8.040,5.882)--(10.849,7.589);
\gppoint{gp mark 2}{(2.422,2.525)}
\gppoint{gp mark 2}{(5.231,4.218)}
\gppoint{gp mark 2}{(8.040,5.882)}
\gppoint{gp mark 2}{(10.849,7.589)}
\gppoint{gp mark 2}{(2.619,7.298)}
\gpcolor{color=gp lt color border}
\node[gp node left] at (3.123,6.836) {$S_{\text{newton}}$};
\gpcolor{rgb color={0.329,0.000,0.800}}
\draw[gp path] (2.391,6.836)--(2.847,6.836);
\draw[gp path] (2.422,2.525)--(5.231,4.235)--(8.040,5.987)--(10.849,7.734);
\gppoint{gp mark 3}{(2.422,2.525)}
\gppoint{gp mark 3}{(5.231,4.235)}
\gppoint{gp mark 3}{(8.040,5.987)}
\gppoint{gp mark 3}{(10.849,7.734)}
\gppoint{gp mark 3}{(2.619,6.836)}
\gpcolor{color=gp lt color border}
\node[gp node left] at (3.123,6.374) {$S_{\text{ideal}}$};
\gpcolor{rgb color={0.800,0.533,0.000}}
\draw[gp path] (2.391,6.374)--(2.847,6.374);
\draw[gp path] (2.422,2.525)--(5.231,4.348)--(8.040,6.171)--(10.849,7.994);
%% coordinates of the plot area
\gpdefrectangularnode{gp plot 1}{\pgfpoint{1.839cm}{1.613cm}}{\pgfpoint{11.671cm}{8.287cm}}
\end{tikzpicture}
%% gnuplot variables
}
    \caption{Computation of shape derivative}
  \end{subfigure}

  \begin{subfigure}[t]{\textwidth}
    \centering
    \begin{tabular}{rrrccc}
    \\
    \toprule
    Procs  & Refs  & DoFs  & Linear solver  & Newton solver   & Linear solver     \\
          &       &        &  (elasticity)      & (shape derivative)    & (shape derivative)\\ 
    \midrule
    40 & 6 & 9,291,330 & 303 & 132 & 544  \\
    160 & 7 & 37,153,922 & 297 & 140 & 567  \\
    640 & 8 & 148,592,898 & 301 & 140 & 625  \\
    2560 & 9 & 594,326,018 & 293 & 140 & 625     \\
     \bottomrule
     \end{tabular}
    \caption{Accumulated iteration counts for the weak scaling study}
    \label{tab:WeakScaling}
  \end{subfigure}

    \caption{Weak scaling: accumulated wallclock time for the first 10 optimization steps (left), speedup relative to 40 cores (right), and accumulated iteration counts (table). Shown are the fine level matrix assembling (\texttt{ass}), multigrid initialization (\texttt{init}) and solution (\texttt{solve}) phase, and the overall newton algorithm (\texttt{newton})}
    \label{fig:WeakScaling}
\end{figure}

\section{Conclusion and outlook}
\label{sec:conclusion}
The results in this article propose a shape optimization algorithm that serves two purposes:
First, towards a modeling of cellular composites, it is possible to maintain the thickness of inter-cellular channel width by promoting shape descent directions with penalized gradients.
We incorporate this as a modification into well-established inner products in a locally linearized description of shape variations.
Second, the gradient penalization prevents early mesh degeneration.
Thereby multigrid scalability can be retained, which is a mandatory ingredient in our simulations towards realistic, large-scale simulations of cell structures in the human skin.
We want to emphasize that our proposed method is not limited to linear elastic models with a limited degree of realism.
The optimization algorithm immediately applies to more complex, potentially nonlinear elasticity models without any changes in methodology.
Under the assumption of hyperelastic material laws, FE models for the human skin are a vivid field of research (e.g., \cite{Querleux2014,Limbert2017MathematicalAC}).
It is thus our next step to apply the proposed optimization algorithm to refined models with experimentally validated material parameters.

\subsection*{Acknowledgements}
This work has been supported by the DFG Priority Program 1648 \emph{Software for Exascale Computing} (SPPEXA) in the project \emph{Exasolvers}. 
We gratefully acknowledge the computing time on the \emph{Statik\&Dynamik}-Cluster at the Department of Civil and Environmental Engineering of the Ruhr University Bochum, Germany.

\bibliographystyle{plain}      % mathematics and physical sciences
\bibliography{citations.bib}   % name your BibTeX data base

\end{document}